\documentclass[11pt]{article}
\usepackage[a4paper,margin=1in]{geometry}
\usepackage{graphicx,xcolor,cite,mathtools}
\usepackage{amssymb,amsmath,amsthm,mathrsfs,paralist,bm,esint,setspace}
\usepackage{tikz}
\usetikzlibrary{arrows.meta,positioning}
\RequirePackage[colorlinks,linkcolor=black,citecolor=blue,urlcolor=blue]{hyperref}
\usepackage{comment}

\newtheorem{theorem}{Theorem}[section]
\newtheorem{lemma}[theorem]{Lemma}
\newtheorem{proposition}[theorem]{Proposition}

\newtheorem{corollary}[theorem]{Corollary}

\theoremstyle{remark}
\newtheorem{remark}[theorem]{Remark}
\theoremstyle{plain}

\newcommand{\Ai}{\operatorname{Ai}}
\newcommand{\Tr}{\operatorname{Tr}}
\newcommand{\Pp}{\mathbb P}
\newcommand{\Ee}{\mathbb E}
\newcommand{\R}{\mathbb R}
\newcommand{\Cov}{\operatorname{Cov}}
\newcommand{\Lip}{\operatorname{Lip}}
\newcommand{\norm}[1]{\left\lVert #1\right\rVert}

\numberwithin{equation}{section}

\title{Gumbel Fluctuations of the Maximum of the Airy$_1$ Process}

\author{Le Chen and Fei Pu}
\date{}

\begin{document}
\maketitle

\begin{abstract}
  We establish the Gumbel limit for the maximum of the Airy$_1$ process over $[0,N]$ as $N\to\infty$, with explicit centering of order $(\log N)^{2/3}$ and scale $(\log N)^{-1/3}$. The proof starts with a two-term expansion for the tail of the maximum over a fixed interval. Its leading term determines the centering, while the second term gives the joint-exceedance asymptotic for touching intervals. We then bound joint exceedances in separated unit intervals, uniformly over gaps up to a fixed multiple of the square root of the threshold. Covariance and mixing estimates control the larger gaps and complete the blocking argument. 
  %The limit can also be deduced from the quantitative persistence estimates of Ferrari and Liu. Our proof is direct and describes the local structure of high exceedances. To the best of our knowledge, the joint-exceedance bound for separated intervals is new.
\end{abstract}

\section{Introduction}

\subsection{The maximum over a long interval}
The Airy$_1$ process $\{\mathcal A_1(x):x\in\mathbb R\}$, derived by Sasamoto \cite{Sas05} in the context of TASEP,  is the universal spatial fluctuation process associated with flat initial data in the Kardar--Parisi--Zhang universality class. It is stationary, has GOE Tracy--Widom one-point marginals up to normalization, and appears as the scaling limit of spatial height profiles in a number of integrable growth and last-passage percolation models; see, for example, \cite{QuastelRemenik2012,QuastelRemenik2013}.
Its one-point distribution is
\[
  \Pp\{\mathcal A_1(0)\le s\}=F_1(2s),
\]
where $F_1$ is the GOE Tracy--Widom distribution function.   In this paper, we study the maximum
\[
  \mathcal M_N:=\max_{0\le x\le N}\mathcal A_1(x)
\]
over a growing interval.
Pu~\cite{Pu2025} proved that, almost surely, $\mathcal M_N/(\log N)^{2/3}$ converges to $\frac12(3/2)^{2/3}$. Basu and Bhattacharjee~\cite{BasuBhattacharjee2024} obtained almost-sure laws for the extrema of both Airy processes. These laws determine the first-order growth of the maximum. A natural question is to describe the fluctuations of the maximum around its leading-order location.
Ferrari and Liu \cite{FL24} studied the related problem of
fixed-level persistence for the Airy$_1$ process and determined
its exponential decay rate as the observation interval grows.

Our main result establishes that the fluctuation of the maximum of the Airy$_1$ process is asymptotically Gumbel.  Throughout the paper, put
\[
  \kappa:=\frac{4\sqrt2}{3},\qquad
  \gamma:=\kappa^{-2/3}=\frac12\left(\frac32\right)^{2/3}.
\]
\begin{theorem}\label{thm:gumbel}
  For $N>1$, define
  \begin{align}\label{eq:aN}
    a_N
    &=\gamma(\log N)^{2/3}
    +\frac{\gamma}{9}\,
    \frac{\log\log N+\log\!\left(\frac{3}{2\pi^3}\right)}{(\log N)^{1/3}},
  \end{align}
  and
  \begin{align}\label{eq:bN}
    b_N=\frac{2\gamma}{3}(\log N)^{-1/3}.
  \end{align}
  Then, for every $z\in\R$,
  \[
    \lim_{N\to\infty}
    \Pp\left\{
    \frac{\mathcal M_N-a_N}{b_N}\le z
    \right\}
    =e^{-e^{-z}}.
  \]
\end{theorem}

Ferrari and Liu~\cite{FL24} do not state Theorem~\ref{thm:gumbel}. Their main theorem concerns a fixed level $c\ge3/2$. However, the error bounds in their proof depend on the level only through a quantity that is bounded away from $1/7$ uniformly in $c$ and tends to zero as $c\to\infty$. Hence, these bounds remain valid when the level grows with the length of the interval, and they also give Theorem~\ref{thm:gumbel}. We sketch this deduction in Remark~\ref{rem:FL}. The proof in this paper is different. It follows the classical extreme-value approach and gives a direct description of high exceedances of $\mathcal A_1$: a sharp tail over every fixed interval and bounds for joint exceedances in touching and separated intervals. We explain our contribution in Section~\ref{SS:contribution}.

A classical counterpart to the present result is the extreme-value theory of stationary Gaussian processes. Let \((X_t)_{t\in \R}\) be a centered stationary Gaussian process with unit variance and correlation function \(R\). Under suitable local regularity and weak-dependence assumptions, for example
\[
R(t)=1-c|t|^\alpha+o(|t|^\alpha)
\qquad\text{as }t\to0,
\]
for some \(c>0\) and \(\alpha\in(0,2]\), together with an appropriate decay condition on \(R(t)\) as \(t\to\infty\), the maximum of \(X\) over a growing interval converges, after centering and scaling, to the Gumbel distribution; see the classical work of Pickands, Berman and Piterbarg \cite{Berman1964,Berman1971,Pickands1969a,Pickands1969b,Piterbarg1996}. The local behavior of \(R\) determines the scale and intensity of high excursion clusters, while the long-range condition ensures their asymptotic independence. The Airy\(_1\) process exhibits a similar interplay between locally Brownian behavior  \cite{QuastelRemenik2012} and decay of dependence \cite{BasuBhattacharjee2024,BasuBusaniFerrari2023, ChenPu2026}, but it is not Gaussian. Consequently, Gaussian comparison inequalities and the standard covariance-based theory do not apply directly: the local extreme intensity and the required anti-clustering estimates must instead be obtained from the Fredholm determinant structure of the Airy$_1$ process.

\subsection{The tail over a fixed interval}
The first ingredient to prove the Gumbel limit in Theorem \ref{thm:gumbel}
is  the sharp bounded-interval tail asymptotic.
For $T>0$ and $M\in\R$, write
\[
  \pi_T(M):=\Pp\left\{\sup_{0\le x\le T}\mathcal A_1(x)>M\right\}.
\]
For $M>0$, set $C_*:=2^{1/4}/\sqrt\pi$ and define
\begin{align}\label{eq:rho}
  \rho(M)
  =
  C_*
  M^{1/4}\exp\{-\kappa M^{3/2}\}.
\end{align}
Maximizing over an interval changes the polynomial factor in the tail. The GOE right-tail asymptotic gives
\[
  p_{\rm pt}(M):=\Pp\{\mathcal A_1(0)>M\}
  \sim\frac{(2M)^{-3/4}}{4\sqrt\pi}e^{-\kappa M^{3/2}};
\]
see~\cite{BaikBuckinghamDiFranco2008}. Comparing this with~\eqref{eq:rho} gives
\[
  \rho(M)\sim8M p_{\rm pt}(M).
\]
Proposition~\ref{lem:sharp-tail} below shows that $\pi_1(M)\sim\rho(M)$ and identifies the second term. Thus, the leading exponential rate agrees with the point tail, but the factor $8M$ shifts the second-order centering. A logarithmic tail estimate alone does not distinguish the two probabilities.

\begin{proposition}\label{lem:sharp-tail}
  For every fixed $T>0$, as $M\to\infty$
  \[
    \pi_T(M)=T\rho(M)+\frac{\rho(M)}{4M}+o\!\left(\frac{\rho(M)}M\right).
  \]
\end{proposition}

\paragraph{The normalization.}
Suppose for a moment that maxima over unit intervals were independent. At the level $u_N(z):=a_N+b_Nz$, a nontrivial limit would require $N\rho(u_N(z))\to e^{-z}$. Taking logarithms suggests
\[
  \kappa u_N(z)^{3/2}
  =\log N+\frac14\log u_N(z)+\log C_*+z+o(1).
\]
At the leading scale $u_N(z)\sim\gamma(\log N)^{2/3}$, the logarithmic prefactor contributes $\frac16\log\log N$. Expanding the left side gives the coefficient $\gamma/9$ in~\eqref{eq:aN}. The slope of $\kappa u^{3/2}$ at $u=\gamma(\log N)^{2/3}$ equals $\frac{3}{2\gamma}(\log N)^{1/3}$, whose reciprocal is $b_N$. Section~\ref{sec:pf} gives the full calculation, including the constant term. The rest of the proof controls dependence between these maxima.

\subsection{Touching and separated blocks}
The second ingredient is the estimate of  joint exceedance probability. 
Let
\[
  Y_j:=\sup_{j\le x\le j+1}\mathcal A_1(x),\qquad j\in\mathbb Z.
\]
Then $\mathcal M_N=\max_{0\le j<N}Y_j$ for integer $N$.

Single-block tails must be supplemented by information about pairs. To see why, consider observations whose exceedance probability at a chosen level is $p_N$, with $Np_N\to\tau>0$. For $N$ independent observations, the no-exceedance probability tends to $e^{-\tau}$. If $N$ is even and each of $N/2$ independent observations is repeated twice, the marginal tail is the same, but
\[
  \Pp\{\text{no exceedance among the repeated observations}\}
  =(1-p_N)^{N/2}\longrightarrow e^{-\tau/2}.
\]
Thus, repeated exceedances can change the limit even when the marginal tail is known.

For the Airy$_1$ block maxima, the second term in Proposition~\ref{lem:sharp-tail} gives the following control of neighboring blocks.

\begin{corollary}\label{lem:touching}
  As $M\to\infty$,
  \[
    \Pp(Y_0>M,Y_1>M)
    \sim\frac{\rho(M)}{4M}.
  \]
\end{corollary}

Indeed, stationarity and inclusion--exclusion give
\[
  \Pp(Y_0>M,Y_1>M)=2\pi_1(M)-\pi_2(M).
\]
The leading terms cancel, leaving $\rho(M)/(4M)$. In particular, $\Pp(Y_1>M\mid Y_0>M)\sim1/(4M)$ tends to zero. Our blocking argument also needs paired exceedances to remain negligible after summing over a growing range of gaps.

For the blocks $[0,1]$ and $[r+1,r+2]$, the gap has length $r$. Our main two-interval estimate is the following.

\begin{proposition}\label{lem:separated}
  Fix  a real number
  $
  \chi\in((1/2)^{1/3},1).
  $
  There exist constants $C,K,M_0>0$, depending only on $\chi$, such that for every $M\ge M_0$ and every integer $r$ with $1\le r\le\chi\sqrt{2M}$,
  \begin{equation}\label{eq:sep-final-bound}
    \Pp(Y_0>M,Y_{r+1}>M)
    \le
    C M^{K}(1+r)^{K}
    \left(e^{-\mathcal J_M(r)}+e^{-2\kappa M^{3/2}}\right),
  \end{equation}
  where
  \begin{equation}\label{eq:JMr}
    \mathcal J_M(r)
    =
    \frac23(2M+r^2)^{3/2}+2Mr+\frac23r^3.
  \end{equation}
  Consequently,
  \begin{align}\label{eq:sum}
    \sum_{j=2}^{\lfloor \chi\sqrt{2M}\rfloor}
    \Pp(Y_0>M,Y_j>M)=o(\rho(M)).
  \end{align}
\end{proposition}

Since $\mathcal J_M(r)\ge\kappa M^{3/2}+2Mr$, Proposition~\ref{lem:separated} gives, for a unit gap and all sufficiently large $M$,
\[
  \Pp(Y_0>M,Y_2>M)\le CM^K\rho(M)e^{-2M},
\]
after enlarging $C,K>0$ if necessary. The second term $e^{-2\kappa M^{3/2}}$ in~\eqref{eq:sep-final-bound} is also covered, since $e^{-\kappa M^{3/2}}\le e^{-2M}$ for large $M$. For touching blocks, the joint probability has the larger order $\rho(M)/M$.

\subsection{From local estimates to the limit}
The proof separates the gaps into ranges with different estimates. With $\chi$ as in Proposition~\ref{lem:separated}, put
\[
  m_M:=\lfloor\chi\sqrt{2M}\rfloor.
\]
Proposition~\ref{lem:separated} controls the sum of joint exceedance probabilities for $2\le j\le m_M$. For the remaining separations of order $\sqrt M$, we discretize each unit interval and use the association of Airy$_1$ proved by Pu~\cite{Pu2025}, together with Newman's inequality~\cite{newman:80:normal}. The covariance estimate of Basu, Busani and Ferrari~\cite{BasuBusaniFerrari2023} gives
\[
  \Cov(\mathcal A_1(0),\mathcal A_1(d))
  =\exp\left\{-\left(\frac43+o(1)\right)d^3\right\}.
\]
The leading covariance exponent at this cutoff satisfies
\[
  \frac43(\chi\sqrt{2M})^3=2\chi^3\kappa M^{3/2}.
\]
Thus, $\chi^3>1/2$ makes the covariance decay faster than $\rho(M)$. Since the two-interval estimate applies for $\chi<1$, both arguments are available near $m_M$.

\begin{figure}[htbp]
  \centering
  \begin{tikzpicture}[x=1cm,y=1cm,>=Latex,font=\small]
    \fill[blue!12] (0,1.4) rectangle (1.8,1.85);
    \fill[blue!12] (6.4,1.4) rectangle (8.2,1.85);
    \draw[thick] (0,1.4)--(8.7,1.4);
    \foreach \x/\lab in {0/0,1.8/1,6.4/j,8.2/{j+1}} {
    \draw (\x,1.32)--(\x,1.5);
    \node[below] at (\x,1.32) {$\lab$};}
    \node at (0.9,1.64) {$Y_0$};
    \node at (7.3,1.64) {$Y_j$};
    \draw[<->] (1.8,2.1)--node[above] {gap $r=j-1$}(6.4,2.1);
    \node[anchor=west] at (9,1.65) {unit intervals};
    \draw[->] (0,-0.05)--(12,-0.05) node[right] {$j$};
    \fill[black!15] (0.2,0) rectangle (1.35,0.35);
    \fill[blue!14] (1.45,0) rectangle (4.8,0.35);
    \fill[teal!18] (4.9,0) rectangle (8.45,0.35);
    \fill[orange!18] (8.55,0) rectangle (11.6,0.35);
    \node at (0.78,0.68) {$j=1$};
    \node at (3.12,0.68) {$2\le j\le m_M$};
    \node at (6.67,0.68) {$m_M<j\le q_N$};
    \node at (10.08,0.68) {$j>q_N$};
    \node[align=center,below] at (0.78,-0.25) {touching\\blocks};
    \node[align=center,below] at (3.12,-0.25) {two-interval\\Fredholm estimate};
    \node[align=center,below] at (6.67,-0.25) {discretization\\and covariance};
    \node[align=center,below] at (10.08,-0.25) {strong\\mixing};
  \end{tikzpicture}
  \caption{Pairs of unit intervals and the estimates used in the proof. The upper diagram shows $j\ge2$ and the gap $r=j-1$. In the lower diagram, $M=u_N(z)$, $m_M=\lfloor\chi\sqrt{2M}\rfloor$, and $A$ is large enough that $q_N>m_M$. The drawing is not to scale; the thresholds indicate where the estimates are applied.}
  \label{Fig:proof-ranges}
\end{figure}
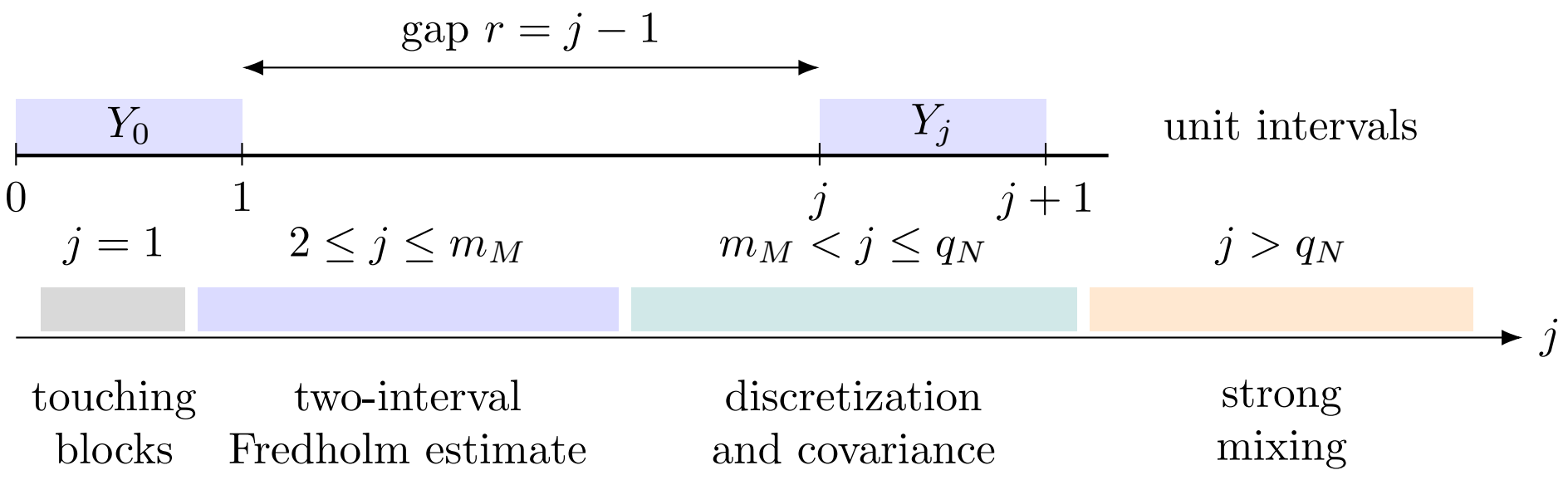

For distant blocks, we use strong mixing: events supported on opposite sides of a gap become nearly independent. Chen and Pu~\cite{ChenPu2026} bound the corresponding $\alpha$-mixing coefficient by $Ce^{-cd^3}$ when the spatial gap has length $d$. We choose
\[
  q_N:=\left\lceil A(\log N)^{1/3}\right\rceil
\]
with $A$ large enough that the accumulated mixing error tends to zero. For fixed $z$ and large $N$, put $M=u_N(z)>0$. Then $q_N$ is of order $\sqrt M$, and we may enlarge $A$ so that $q_N>m_M$. Figure~\ref{Fig:proof-ranges} summarizes the ranges.

Combining the local estimates gives
\[
  N\sum_{j=1}^{q_N}\Pp(Y_0>u_N(z),Y_j>u_N(z))\longrightarrow0.
\]
We then divide the observation interval into long blocks separated by gaps of length $q_N$. The local bound controls repeated exceedances inside a long block, while strong mixing controls dependence between different long blocks. The discarded gaps have negligible exceedance probability. We take a liminf and a limsup in $N$ with the number of long blocks fixed. Letting this number increase makes the two bounds for $\Pp\{\mathcal M_N\le u_N(z)\}$ agree and proves Theorem~\ref{thm:gumbel}.

\subsection{Contribution and further questions}\label{SS:contribution}
Let us summarize what is new here and what is due to Ferrari and Liu~\cite{FL24}. Their quantitative persistence estimates give the Gumbel limit in Theorem~\ref{thm:gumbel}, as explained in Remark~\ref{rem:FL}. For $T\ge1$, the two-term expansion in Proposition~\ref{lem:sharp-tail}, with an absolute remainder bound, can also be extracted from intermediate estimates in their proof, although not from their stated results. Moreover, the first trace identity in our proof of Proposition~\ref{lem:sharp-tail} is a fixed-interval counterpart of~\cite[Proposition 2.5]{FL24}.

The contributions of this paper are the following.
\begin{enumerate}[(i)]
  \item A direct proof of the Gumbel limit through the local structure of high exceedances, in the spirit of the extreme-value theory for stationary Gaussian processes. The proof does not use the bounds of~\cite{FL24} for higher powers of the persistence kernel.
  \item The two-term tail expansion in Proposition~\ref{lem:sharp-tail} for every fixed $T>0$, with an explicit signed remainder, and the resulting asymptotic for touching blocks in Corollary~\ref{lem:touching}.
  \item The joint-exceedance bound for two separated intervals in Proposition~\ref{lem:separated}. It is uniform over gaps up to $\chi\sqrt{2M}$. The paper~\cite{FL24} treats a single interval and does not contain such an estimate.
  \item The control of all remaining gaps. It combines Proposition~\ref{lem:separated} with association, the covariance decay in~\cite{BasuBusaniFerrari2023} and the mixing bound in~\cite{ChenPu2026}.
\end{enumerate}
Together, (ii)--(iv) show that high exceedances in different unit intervals do not cluster at the extreme-value scale. Indeed, $\Pp(Y_1>M\mid Y_0>M)\sim1/(4M)$, and the separated pairs satisfy~\eqref{eq:sum}.

Our separated-block estimate is an upper bound. It does not determine a matching asymptotic throughout the range of growing gaps. A further question is the joint local profile of the process near several extreme peaks, which maxima over unit intervals do not capture.

The paper is organized as follows. Section~\ref{sec:tail} proves the fixed-interval expansion and the touching-block asymptotic. Section~\ref{sec:short} establishes the two-interval estimate through a Fredholm determinant expansion. Section~\ref{sec:mes} treats the mesoscopic range. Section~\ref{sec:pf} completes the normalization and blocking argument.  The appendices supply the continuum determinant limit and the auxiliary operator estimates.

%\bigskip
%\setcounter{tocdepth}{2}
%\begingroup\small\tableofcontents\endgroup
%\clearpage

\section{Sharp tail over a bounded interval} \label{sec:tail}

We now prove Proposition~\ref{lem:sharp-tail} and Corollary~\ref{lem:touching}. The first step is to specify the Fredholm determinant and the weighted kernels used below.

\subsection{Fredholm conventions}\label{SS:fredholm-conventions}
Throughout the proofs, we use the multiplication operators $U_a$ ($a>0$) and $\Phi$ defined by
\begin{align}\label{eq:UPhi}
  \qquad U_af(x)=e^{-ax}f(x), \qquad \text{and} \qquad
  \Phi f(x)=\phi(x)f(x),\quad \text{with}\,\, \, \phi(x)=(1+x^2)^{1/2}.
\end{align}
The factor $\Phi$ supplies a polynomial weight in the factorizations below. The operators $U_a$ and $U_a^{-1}$ are unbounded on $L^2(\R)$. Products involving them denote weighted kernels, defined first on compactly supported functions and then extended using the norm bounds in Appendix~\ref{S:operator-appendix}. We call an exponential weight admissible for a kernel if the weighted kernel extends to a trace-class operator.

Suppose $K=VW$, where $V$ and $W$ are Hilbert--Schmidt. We use the kernel representative
\[
  K(x,y)=\int_\R V(x,z)W(z,y)\,dz.
\]
The trace formula on exterior powers gives
\[
  \Tr\left(\wedge^n(VW)\right)
  =\frac1{n!}\int_{\R^n}\det[K(x_i,x_j)]_{i,j=1}^n\,d\mathbf x;
\]
see~\cite[Chapter 3]{Simon2005}. By the continuous Cauchy--Binet identity,
\[
  \det[K(x_i,x_j)]_{i,j=1}^n
  =\frac1{n!}\int_{\R^n}
  \det[V(x_i,z_j)]_{i,j=1}^n
  \det[W(z_i,x_j)]_{i,j=1}^n\,d\mathbf z.
\]
Both determinants on the right are square-integrable in the $2n$ variables. Thus, Cauchy--Schwarz shows that the iterated integral converges absolutely. For two such admissible exponential weights, the row and column factors cancel in each kernel determinant. Hence, the Fredholm coefficients and the determinant are independent of the weight.

\subsection{The fixed-interval expansion}
\begin{proof}[Proof of Proposition~\ref{lem:sharp-tail}] Throughout the proof $T>0$ is fixed. Take $M$ sufficiently large that
  $
  a:=\sqrt{2M}\ge T+1.
  $
  Recall that $P_t=e^{t\Delta}$ denotes the heat semigroup, with kernel
  \[
    p_t(x-z)=\frac1{\sqrt{4\pi t}}\exp\left\{-\frac{(x-z)^2}{4t}\right\},\qquad t>0.
  \]

  \medskip
  \noindent\textbf{Step~1}.
  The continuum-statistics formula of~\cite[Theorem 4]{QuastelRemenik2012} gives
  \begin{equation}\label{eq:QR-oneblock}
    \Pp\{\mathcal A_1(x)\le M,\ 0\le x\le T\}
    =
    \det\left(I-B_0+\Lambda^M_{[0,T]}P_{-T}B_0\right),
  \end{equation}
  where $B_0$ has kernel $B_0(x,y)=\Ai(x+y)$ and the kernel $\Lambda^M_{[0,T]}$ has the form
  \begin{align*}
    p_T(x-z) \Pp_{x,z}^T\{b(t)\leq M, t\in [0, T]\}.
  \end{align*}
  Here, $b$ is a Brownian bridge starting from $x$ at time $0$ and ending at $z$ at time $T$ with diffusion coefficient $2$.
  This probability equals $1- e^{-(x-M)(z-M)/T}$
  if $x\leq M, z\leq M$ and $0$ otherwise (see~\cite[Chapter IV, No.~26, p.~67]{BoS02}).
  Hence, we have
  \begin{align}\label{eq:Lambda}
    \Lambda^M_{[0,T]}(x,z)
    =\mathbf 1_{\{x\le M,z\le M\}}
    \left[p_T(x-z)-p_T(2M-x-z)\right].
  \end{align}
  For $T>0$, we shall use the explicit identity
  \begin{equation}\label{eq:backward-airy-kernel}
    (P_{-T}B_0)(z,y)
    =e^{-2T^3/3-T(z+y)}\Ai(z+y+T^2)
  \end{equation}
  (see~\cite[(1.9)]{QuastelRemenik2012}).
  Set
  $
  Q_{M,T}:=P_T-\Lambda^M_{[0,T]}
  $. Then
  \begin{equation}\label{eq:Q-piecewise-L21}
    Q_{M,T}(x,z)
    =\begin{cases}
      \displaystyle p_T(2M-x-z),&x\le M,\ z\le M,\\[1mm]
      \displaystyle p_T(x-z),&\text{at least one of $x,z$ is $>M$}.
    \end{cases}
  \end{equation}

  Denote $K_{M,T}:=Q_{M,T}P_{-T}B_0$.
  By the scalar semigroup identity in~\cite[Proposition 1.2]{QuastelRemenik2012},
  $P_TP_{-T}B_0=B_0$ pointwise, and hence
  $K_{M,T}=B_0-\Lambda^M_{[0,T]}P_{-T}B_0$.
  These scalar compositions are absolutely convergent because $0\le\Lambda^M_{[0,T]}\le P_T$ and the Gaussian factor dominates the exponential bound on the backward Airy kernel.
  We realize the Fredholm determinant of this kernel through a trace-class weighted operator.  Let $U_a$ and $\phi$ be defined as in~\eqref{eq:UPhi}.
  Introduce the two kernels
  \begin{align}
    V_{T,a}(x,z)
    &:=e^{-ax}Q_{M,T}(x,z)e^{-az}\phi(z),\label{eq:V-L21}\\
    W_{T,a}(z,y)
    &:=\phi(z)^{-1}e^{az}(P_{-T}B_0)(z,y)e^{ay}.\label{eq:W-L21}
  \end{align}
  By Lemmas~\ref{lem:V} and~\ref{lem:W}, the kernels $V_{T,a}$ and $W_{T,a}$ are Hilbert--Schmidt under the imposed condition $a\ge T+1$, and consequently
  \[
    \widetilde A_{M,T}:=V_{T,a}W_{T,a}
    \] is trace class.
    Remark 1.1 of~\cite{QuastelRemenik2012} also interprets~\eqref{eq:QR-oneblock} as the Fredholm expansion of this kernel. By the coefficient formula above, that series equals $\det(I-\widetilde A_{M,T})$. Thus, \begin{equation}\label{eq:p-det-L21}
      \pi_T(M)=1-\det(I-\widetilde A_{M,T}).
    \end{equation}

    \medskip
    \noindent\textbf{Step~2}.
    Since $V_{T,a}$ and $W_{T,a}$ are Hilbert--Schmidt,
    and the standard trace formula for a product of two Hilbert--Schmidt operators gives
    \begin{align}
      \Tr\widetilde A_{M,T}
      &=\iint_{\R^2}V_{T,a}(x,z)W_{T,a}(z,x)\,dz\,dx
      =\iint_{\R^2}Q_{M,T}(x,z)(P_{-T}B_0)(z,x)\,dz\,dx.\label{eq:trace-QW-L21}
    \end{align}
    We use~\eqref{eq:backward-airy-kernel} and split~\eqref{eq:trace-QW-L21} according to the two regions in~\eqref{eq:Q-piecewise-L21}:
    \begin{equation*}\label{eq:trace-split-L21}
      \Tr\widetilde A_{M,T}=I_{M,T}+R_{M,T},
    \end{equation*}
    where
    \begin{align*}
      I_{M,T}&=  \frac{e^{-2T^3/3}}{\sqrt{4\pi T}}\int_{-\infty}^M \int_{-\infty}^M
      \exp\left\{-\frac{(2M-x-z)^2}{4T}-T(x+z)\right\}
      \Ai(x+z+T^2) dxdz,\\
      R_{M,T}&=   \frac{e^{-2T^3/3}}{\sqrt{4\pi T}}\iint_{x>M \, \text{or}\, z>M}
      \exp\left\{-\frac{(x-z)^2}{4T}-T(x+z)\right\}
      \Ai(x+z+T^2) dxdz.
    \end{align*}
    A change of variable $s=x+z$ yields
    \begin{equation}\label{eq:IMT}
      I_{M,T}
      =
      \frac{e^{-2T^3/3}}{2\sqrt{\pi T}}
      \int_{-\infty}^{2M}
      (2M-s)e^{-(2M-s)^2/(4T)-Ts}\Ai(s+T^2)\,ds.
    \end{equation}
    We also simplify the expression of $R_{M,T}$.  Introduce
    \[
      s=x+z,\qquad d=x-z,
      \qquad dx\,dz=\frac12\,ds\,dd.
    \]
    If $s\le2M$ and $b:=2M-s\ge0$, then the condition that at least one endpoint is above $M$ is exactly $|d|>b$.  If $s>2M$, every $d\in\R$ belongs to the complementary region.  Therefore, \begin{align*}
      R_{M,T}=R_{M,T}^{(-)}+R_{M,T}^{(+)},
    \end{align*}
    where
    \begin{align}
      R_{M,T}^{(-)}
      &=\frac{e^{-2T^3/3}}{2\sqrt{\pi T}}
      \int_{-\infty}^{2M}e^{-Ts}\Ai(s+T^2)
      \left(\int_{2M-s}^{\infty}e^{-d^2/(4T)}\,dd\right)ds,
      \label{eq:Rminus-exact-L21}\\
      R_{M,T}^{(+)}
      &=\frac{e^{-2T^3/3}}2
      \int_{2M}^{\infty}e^{-Ts}\Ai(s+T^2)\,ds.
      \label{eq:Rplus-exact-L21}
    \end{align}
    The factor in~\eqref{eq:Rminus-exact-L21} comes from the two symmetric tails $d>b$ and $d<-b$, which cancel the Jacobian factor $1/2$.

    Put $A=2M$ and $g_T(s)=e^{-2T^3/3-Ts}\Ai(s+T^2)$, and write
    $\overline P_T(v)=\int_v^\infty p_T(d)\,dd$ for the Gaussian tail, for every real $v$.
    The exact split above reads
    \begin{align*}
      I_{M,T}&=\int_{-\infty}^A(A-s)p_T(A-s)g_T(s)\,ds,\\
      R_{M,T}^{(-)}&=\int_{-\infty}^A\overline P_T(A-s)g_T(s)\,ds,\\
      R_{M,T}^{(+)}&=\frac12\int_A^\infty g_T(s)\,ds.
    \end{align*}
    The scalar heat--Airy identity of Quastel--Remenik~\cite[Proposition 1.2]{QuastelRemenik2012} is
    \[
      \int_\R p_T(x-s)g_T(s)\,ds=\Ai(x).
    \]
    This scalar convolution is absolutely convergent. Since $|g_T(s)|\le C_Te^{-Ts}$ and Gaussian derivatives dominate this exponential on every compact set of $x$, differentiation under the integral gives
    \[
      \int_\R(A-s)p_T(A-s)g_T(s)\,ds=-2T\Ai'(A).
    \]
    For the second full-line identity, use $\overline P_T(A-s)=\int_\R p_T(d)\mathbf1_{\{s+d\ge A\}}\,dd$. The double integral is absolutely convergent because
    \begin{align*}
      \int_\R p_T(d)\int_{A-d}^\infty |g_T(s)|\,ds\,dd
      &\le \frac{C}{T}e^{-TA-2T^3/3}\int_\R p_T(d)e^{Td}\,dd\\
      &=\frac{C}{T}e^{-TA+T^3/3}<\infty.
    \end{align*}
    Fubini, followed by $u=s+d$ and the heat--Airy identity, therefore yields
    \[
      \int_\R\overline P_T(A-s)g_T(s)\,ds
      =\int_A^\infty\int_\R p_T(d)g_T(u-d)\,dd\,du
      =\int_A^\infty\Ai(u)\,du.
    \]
    Subtract the parts with $s>A$ from these two identities and add $R_{M,T}^{(+)}$. Evenness of the heat kernel gives
    \begin{equation}\label{eq:trace-heat-exact}
      \Tr\widetilde A_{M,T}
      =-2T\Ai'(A)+\int_A^\infty\Ai(u)\,du+\mathcal E_{M,T},
    \end{equation}
    where
    \[
      \mathcal E_{M,T}=\int_0^\infty
      \left(bp_T(b)-\int_0^b p_T(d)\,dd\right)g_T(A+b)\,db.
    \]
    For $b\ge0$ the heat kernel is decreasing, so the expression in parentheses lies in $[-1/2,0]$. Since $g_T(A+b)>0$ for $A>0$, we have
    \[
      -\frac12\int_A^\infty g_T(s)\,ds\le\mathcal E_{M,T}\le0.
    \]
    For $A=a^2\ge1$ and $b\geq0$,  the  positive-argument Airy bound $0<\Ai(x) \leq C x^{-1/4}e^{-\frac23 x^{3/2}}$ for $x>0$ gives
    \begin{align*}
    \Ai(a^2+b+T^2) \leq C a^{-1/2} e^{-\frac23(a^2+b+T^2)^{3/2}}.
    \end{align*}
    Because the function $x\to x^{3/2}$ is convex on $[0, \infty)$,  its tangent line at \(x=a^2\) lies below the graph and hence we have
    \begin{align*}
    (a^2+b+T^2)^{3/2} \geq a^3+ \frac32 a(b+T^2) \geq a^3 +\frac32 ab. 
    \end{align*}
It follows that    
\begin{align*}
      g_T(A+b)
      &\le C_Ta^{-1/2}e^{-TA-2a^3/3}e^{-(T+a)b},\qquad b\ge0.
          \end{align*}
    As a consequence, 
    \begin{align*}
      |\mathcal E_{M,T}|
      &\le C_T\frac{a^{-1/2}}{a+T}e^{-TA-2a^3/3}
      =O_T\!\left(a^{-3/2}e^{-Ta^2-2a^3/3}\right)
      =o\!\left(\frac{\rho(M)}M\right).
    \end{align*}
    
    Finally, the classical Airy and Airy-derivative asymptotics (see~\cite[Section 10.4]{AS64}) give
    \[
      -\Ai'(A)=\frac{A^{1/4}e^{-2A^{3/2}/3}}{2\sqrt\pi}
      \bigl(1+O(A^{-3/2})\bigr),\qquad
      \int_A^\infty\Ai(u)\,du
      =\frac{A^{-3/4}e^{-2A^{3/2}/3}}{2\sqrt\pi}\bigl(1+o(1)\bigr).
    \]
    The second relation follows, for example, by l'Hopital's rule from the usual Airy asymptotic. Thus, \[
    -2T\Ai'(2M)=T\rho(M)+O_T(M^{-3/2}\rho(M)),\qquad
    \int_{2M}^\infty\Ai(u)\,du=\frac{\rho(M)}{4M}(1+o(1)).
  \]
  Together with~\eqref{eq:trace-heat-exact}, this proves
  \begin{equation}\label{eq:traceM}
    \Tr\widetilde A_{M,T}=T\rho(M)+\frac{\rho(M)}{4M}
    +o\!\left(\frac{\rho(M)}M\right).
  \end{equation}
  The identity~\eqref{eq:trace-heat-exact} is a fixed-interval counterpart of the first trace identity in Ferrari--Liu~\cite[Proposition 2.5]{FL24}. The fixed-$T$, large-$M$ remainder estimate above is the one needed here.

  \medskip
  \noindent\textbf{Step~3}.
  Lemmas~\ref{lem:W} and~\ref{lem:V}  prove the quantitative trace-norm estimate
  \begin{equation}\label{eq:trace-norm-L21}
    \|\widetilde A_{M,T}\|_1 \leq \|V_{T, a}\|_2 \|W_{T, a}\|_2
    \le C_TM^{K}e^{-2a^3/3}
  \end{equation}
  for some constant $K>0$. The exponent is
  $(Ta^2-a^3)+(a^3/3-Ta^2)=-2a^3/3$, so the $T$-dependent exponential factors cancel. This yields
  \begin{equation}\label{eq:trace-norm-square-L21}
    \|\widetilde A_{M,T}\|_1^2
    =o\left(\frac{\rho(M)}M\right).
  \end{equation}
  For a trace-class operator $A$, we have
  \begin{align}\label{eq:tracein}
    \left|1-\det(I-A)-\Tr A\right|
    \le\frac12\|A\|_1^2e^{\|A\|_1}
  \end{align}
  (see~\cite[Lemma 3.3]{Simon2005}).
  Applying this with $A=\widetilde A_{M,T}$ and using~\eqref{eq:trace-norm-square-L21},
  \begin{equation}\label{eq:det-remainder-L21}
    1-\det(I-\widetilde A_{M,T})
    =\Tr\widetilde A_{M,T}
    +o\left(\frac{\rho(M)}M\right).
  \end{equation}
  Finally, we combine~\eqref{eq:p-det-L21}, \eqref{eq:det-remainder-L21} with~\eqref{eq:traceM} to complete the proof of Proposition~\ref{lem:sharp-tail}.
  This completes the proof of Proposition~\ref{lem:sharp-tail}.
\end{proof}

\subsection{Touching intervals}
We apply the fixed-interval expansion at lengths one and two.

\begin{proof}[Proof of Corollary~\ref{lem:touching}]
  By stationarity and inclusion--exclusion,
  \[
    \Pp(Y_0>M,Y_1>M)=2\pi_1(M)-\pi_2(M).
  \]
  Apply Proposition~\ref{lem:sharp-tail} with $T=1$ and $T=2$.
  This completes the proof of Corollary~\ref{lem:touching}.
\end{proof}

\section{Two-interval estimates and short-range anti-clustering}\label{sec:short}

\subsection{The two-interval determinant}\label{SS:two-interval}
We now prove Proposition~\ref{lem:separated}. We start with a determinant formula for the joint exceedance probability
$\mathbb P(Y_0>M,Y_j>M)$ in the stated range of gaps. The adjacent case $j=1$ was already treated in Corollary~\ref{lem:touching} by means of the sharp bounded-interval expansion. We first give a Fredholm determinant formula for the two-interval hitting probability of the Airy$_1$ process.  Recall $\Lambda_{[0, T]}^M$ given in~\eqref{eq:Lambda}.

\begin{lemma}\label{lem:two-interval-determinant}
  Fix $T>0$, $r>0$, and $M>0$ with $\sqrt{2M}\ge 2T+r+1$, and put
  $
  S:=2T+r$, $\Lambda:=\Lambda^M_{[0,T]}.
  $
  Then
  \begin{align}
    &\Pp\left\{\mathcal A_1(x)\le M\text{ for every }
    x\in[0,T]\cup[T+r,2T+r]\right\}\notag\\
    &\qquad=
    \det\left(I-B_0+\Lambda P_r\Lambda P_{-S}B_0\right).
    \label{eq:two-interval-continuum}
  \end{align}
\end{lemma}

\begin{remark}
  The determinant in~\eqref{eq:two-interval-continuum} is understood through a trace-class conjugate of the operator inside the Fredholm expansion.  Equivalently, define
  \begin{equation}\label{eq:Kcup-def-main}
    K_\cup:=\bigl(P_TP_rP_T-\Lambda P_r\Lambda\bigr)P_{-S}B_0.
  \end{equation}
  For $U=U_a$ with $a=\sqrt{2M}\ge 2T+r+1$, the kernel
  $\widetilde K_\cup:=UK_\cup U^{-1}$ is trace class, as proved in the appendix, and one has
  \begin{equation}\label{eq:two-interval-Kcup-det}
    \Pp\left\{\mathcal A_1(x)\le M\text{ on }
    [0,T]\cup[T+r,2T+r]\right\}
    =\det(I-\widetilde K_\cup).
  \end{equation}
  The proof follows the discretization argument of~\cite[Theorem~4]{QuastelRemenik2012}, together with the finite-dimensional formula in~\cite[Theorem~1]{QuastelRemenik2012}.
\end{remark}

\subsection{Joint exceedances and the interaction operator}\label{SS:joint}
In light of Lemma~\ref{lem:two-interval-determinant}, we will derive a formula for
the joint exceedance probability $\Pp\{Y_0>M, Y_{r+1}>M\}$ by inclusion-exclusion.  Set
\begin{align*}
  \Lambda=\Lambda_{[0,1]}^M, \quad \text{and}\quad  Q=P_1-\Lambda.
\end{align*}
Define
\begin{equation}\label{eq:K12}
  K_1:=QP_{-1}B_0, \qquad K_2:=P_{1+r}QP_{-(2+r)}B_0.
\end{equation}
Note that the factor $P_{1+r}$ in $K_2$ is the unconstrained forward propagation from time $0$ to the beginning of the second block, while $P_{-(2+r)}B_0$ is the closing factor from the final time $2+r$ back to time $0$.

\begin{lemma}\label{lem:joint}
  Fix $M>0$ and $r\in \mathbb{N}$ with $r\ge1$ and $\sqrt{2M}\ge r+3$. Let $a=\sqrt{2M}$.
  Then we have
  \begin{equation}\label{eq:joint-det-inclusion}
    \Pp\{Y_0>M, Y_{r+1}>M\}
    =1-\det(I-\widetilde K_1)-\det(I-\widetilde K_2)
    +\det(I-\widetilde K_\cup),
  \end{equation}
  where $\tilde{K}_1, \tilde{K}_2, \tilde{K}_\cup$ are the trace conjugates with respect to $U_a$ of the operators $K_1, K_2, K_\cup$.  In addition, we have
  \begin{equation}\label{eq:Kcup-conj-decomp-detail}
    \widetilde K_\cup
    =\widetilde K_1+\widetilde K_2-\widetilde C_r
  \end{equation}
  with $C_r=QP_rQP_{-(2+r)}B_0$.

\end{lemma}
\begin{proof}
  Denote
  $
  A=\{Y_0>M\}$ and $B=\{Y_{r+1}>M\}.
  $ Formula (1.14) of~\cite{QuastelRemenik2012}, together with time-translation invariance  gives
  \[
    \Pp(B^c)=\det\!\left(I-QP_{-1}B_0\right),
  \]
  where, as throughout, this means the Fredholm determinant represented by the appropriate trace-class conjugate. Thus, the direct one-block operator is $QP_{-1}B_0$. The operator $K_2$ in~\eqref{eq:K12} is instead the representation naturally adapted to the common interval $[0,2+r]$. We now verify  that the two representations have the same Fredholm determinant.

  Recall the operators $U_a$ and $\Phi$ introduced in~\eqref{eq:UPhi}.
  Let
  \[
    L:=1+r,\qquad S:=2+r=L+1.
  \]
  Define
  \[
    A_{r,a}:=(U_aP_LU_a^{-1})(U_aQU_a\Phi),\qquad
    B_{r,a}:=\Phi^{-1}U_a^{-1}P_{-S}B_0U_a^{-1}.
  \]
  $U_aP_LU_a^{-1}$ is bounded by Lemma~\ref{lem:2to2}.  Moreover, $U_aQU_a\Phi$ and $B_{r,a}$ are Hilbert--Schmidt by Lemmas~\ref{lem:W} and~\ref{lem:V}.  Hence, $A_{r,a}$ is Hilbert--Schmidt. In addition,
  \[
    A_{r,a}B_{r,a}
    =U_aP_LQP_{-S}B_0U_a^{-1}
    =U_aK_2U_a^{-1}=:\widetilde K_2.
  \]
  We  use the following identity for Fredholm determinants,
  \begin{align}\label{eq:commute}
    \det(I-AB)=\det(I-BA),
  \end{align}
  valid for Hilbert--Schmidt $A,B$; see Simon~\cite[Chapter~3]{Simon2005}.  Hence, \[
  \det(I-\widetilde K_2)=\det(I-B_{r,a}A_{r,a}).
\]

We next identify the product on the right.  Using the scalar heat--Airy identity $P_{-S}B_0P_L=P_{-1}B_0$ from~\cite[Proposition 1.2]{QuastelRemenik2012}, with $S-L=1$, we obtain
\[
  B_{r,a}A_{r,a}
  =\Phi^{-1}U_a^{-1}P_{-1}B_0Q\,U_a\Phi.
\]
Now define
\[
  C_a:=U_aQU_a\Phi,\qquad
  D_a:=\Phi^{-1}U_a^{-1}P_{-1}B_0U_a^{-1}.
\]
Then $C_a$ and $D_a$ are Hilbert--Schmidt by Lemmas~\ref{lem:W} and~\ref{lem:V}. Moreover, we have
\[
  C_aD_a=U_aQP_{-1}B_0U_a^{-1}=:\widetilde K_1,
\]
while
\[
  D_aC_a=\Phi^{-1}U_a^{-1}P_{-1}B_0Q\,U_a\Phi=B_{r,a}A_{r,a}.
\]
Applying the  identity~\eqref{eq:commute} once more,
\[
  \det(I-B_{r,a}A_{r,a})=\det(I-C_aD_a)=\det(I-\widetilde K_1).
\]
Therefore, \[
\det(I-\widetilde K_2)=\det(I-\widetilde K_1),
\]
with $\widetilde K_1=U_aK_1U_a^{-1}$ and $\widetilde K_2=U_aK_2U_a^{-1}$.  Thus, we have
\[
  \Pp(A^c)=\det(I-\widetilde K_1),\qquad
  \Pp(B^c)=\det(I-\widetilde K_2).
\]
By Lemma~\ref{lem:two-interval-determinant}, with $T=1$ and $S=2+r$,
\[
  \Pp(A^c\cap B^c)=\det(I-\widetilde K_\cup).
\]
Since
$
\Pp(A\cap B)=1-\Pp(A^c)-\Pp(B^c)+\Pp(A^c\cap B^c),
$
we obtain the identity~\eqref{eq:joint-det-inclusion}.

We next identify the interaction term in $K_\cup$. By~\eqref{eq:Kcup-def-main},
\[
  K_\cup=\bigl(P_1P_rP_1-\Lambda P_r\Lambda\bigr)P_{-(2+r)}B_0.
\]
Using $\Lambda=P_1-Q$,
\begin{align*}
  P_1P_rP_1-\Lambda P_r\Lambda
  &=P_1P_rP_1-(P_1-Q)P_r(P_1-Q)\\
  &=QP_rP_1+P_1P_rQ-QP_rQ.
\end{align*}
Multiplying on the right by $P_{-(2+r)}B_0$ and using the semigroup identities
\[
  P_rP_1P_{-(2+r)}B_0=P_{-1}B_0,
  \qquad
  P_1P_r=P_{1+r},
\]
gives
\[
  K_\cup=K_1+K_2-C_r,
\]
where, consistently with~\eqref{eq:K12},
\[
  C_r=QP_rQP_{-(2+r)}B_0.
\]
Thus, $C_r$ is exactly the interaction term between the two constrained blocks.

Conjugating the remaining kernels by the same multiplication operator $U_a$ gives~\eqref{eq:Kcup-conj-decomp-detail}.
This completes the proof of Lemma~\ref{lem:joint}.
\end{proof}

\subsection{Proof of the joint-exceedance bound}\label{SS:joint-bound}
We now prove Proposition~\ref{lem:separated}. The determinant formula expresses the joint probability as a trace term plus a remainder. We first use the bounds for these two terms to prove~\eqref{eq:sep-final-bound} and then sum over gaps to obtain~\eqref{eq:sum}. Subsections~\ref{SS:trace-sectors}--\ref{SS:trace-norms} prove the required estimates.

\begin{proof}[Proof of Proposition~\ref{lem:separated}]
  Enlarge $M_0$ so that $a=\sqrt{2M}\ge 3/(1-\chi)$. Then $a-(r+2)\ge1$ throughout the stated range, so all weighted factorizations below are admissible. We also take $M_0$ larger than the thresholds in the auxiliary lemmas proved below.

  \noindent\textbf{Step~1}.
  Denote
  $
  A=\{Y_0>M\}$ and $B=\{Y_{r+1}>M\}.
  $
  For a trace-class operator $H$ we use~\eqref{eq:tracein} to write
  \begin{equation}\label{eq:det-linear-remainder-detail}
    \det(I-H)=1-\Tr H+R(H),
  \end{equation}
  with the standard bound
  \begin{equation}\label{eq:det-rem-bound-detail}
    |R(H)|
    \le \frac12\|H\|_1^2 e^{\|H\|_1}.
  \end{equation}
  Applying~\eqref{eq:det-linear-remainder-detail} to the three determinants in~\eqref{eq:joint-det-inclusion}, we obtain
  \begin{align*}
    \Pp(A\cap B)
    &=1-\{1-\Tr\widetilde K_1+R(\widetilde K_1)\}
    -\{1-\Tr\widetilde K_2+R(\widetilde K_2)\}\\
    &\quad+\{1-\Tr\widetilde K_\cup+R(\widetilde K_\cup)\}\\
    &=\Tr\widetilde K_1+\Tr\widetilde K_2-\Tr\widetilde K_\cup
    -R(\widetilde K_1)-R(\widetilde K_2)+R(\widetilde K_\cup).
  \end{align*}
  By~\eqref{eq:Kcup-conj-decomp-detail},
  \[
    \Tr\widetilde K_1+\Tr\widetilde K_2-\Tr\widetilde K_\cup
    =\Tr\widetilde C_r.
  \]
  Hence, \begin{equation}\label{eq:joint-trace-plus-rem-detail}
    \Pp(A\cap B)
    =\Tr\widetilde C_r+\mathcal R_{M,r},
  \end{equation}
  where
  \[
    \mathcal R_{M,r}
    :=-R(\widetilde K_1)-R(\widetilde K_2)+R(\widetilde K_\cup).
  \]
  The one-block traces cancel. It remains to estimate the interaction trace and the determinant remainder.

  The sector estimates below give~\eqref{eq:traceC-main-direct}, and Lemma~\ref{L:joint-trace-norms} gives the required trace-norm bound. Both hold uniformly over the range in the proposition. We now combine these estimates.
  \medskip
  \noindent\textbf{Step~2}. We complete the Fredholm estimate.
  Recall $\frac23a^3=\kappa M^{3/2}$.
  Set
  \[
    \eta_{M,r}:=C M^K(1+r)^K e^{-\kappa M^{3/2}}.
  \]
  Because $r\le\chi a=O(\sqrt M)$, the polynomial factor is at most polynomial in $M$, and therefore $\eta_{M,r}\to0$ uniformly in the stated range.  For large $M$ we may assume $\eta_{M,r}\le1$.  From~\eqref{eq:det-rem-bound-detail} and~\eqref{eq:crude-norm-main-direct},
  \begin{align*}
    |\mathcal R_{M,r}|
    &\le \frac12\Bigl(
    \|\widetilde K_1\|_1^2e^{\|\widetilde K_1\|_1}
    +\|\widetilde K_2\|_1^2e^{\|\widetilde K_2\|_1}
    +\|\widetilde K_\cup\|_1^2e^{\|\widetilde K_\cup\|_1}
    \Bigr)\\
    &\le C M^{K'}(1+r)^{K'}e^{-2\kappa M^{3/2}}.
  \end{align*}
  Combining this with~\eqref{eq:joint-trace-plus-rem-detail} and the direct  trace estimate~\eqref{eq:traceC-main-direct}, and enlarging $K$ if necessary, gives
  \[
    \Pp(Y_0>M,Y_{r+1}>M)
    \le C M^K(1+r)^K
    \left(e^{-\mathcal J_M(r)}+e^{-2\kappa M^{3/2}}\right),
  \]
  which proves~\eqref{eq:sep-final-bound}.

  \medskip
  \noindent\textbf{Step~3}. We finish the proof of~\eqref{eq:sum}.
  First,
  since
  $
  \kappa M^{3/2}=\frac23(2M)^{3/2},
  $
  we have
  \begin{align*}
    \mathcal J_M(r)-\kappa M^{3/2}
    &=\frac23\Bigl[(2M+r^2)^{3/2}-(2M)^{3/2}\Bigr]
    +2Mr+\frac23r^3
    \ge 2Mr,
  \end{align*}
  because the two omitted terms are nonnegative. Therefore, \begin{equation}\label{eq:exp-J-rho-detail}
    e^{-\mathcal J_M(r)}
    \le e^{-\kappa M^{3/2}}e^{-2Mr}.
  \end{equation}

  Recall
  $
  \rho(M)=C_*M^{1/4}e^{-\kappa M^{3/2}}
  $ as in~\eqref{eq:rho}.
  Dividing the first term in~\eqref{eq:sep-final-bound} by $\rho(M)$ and using~\eqref{eq:exp-J-rho-detail} gives, uniformly in $r\ge1$,
  \begin{equation}\label{eq:first-sum-ratio-detail}
    \frac{C M^K(1+r)^K e^{-\mathcal J_M(r)}}{\rho(M)}
    \le C M^{K-1/4}(1+r)^K e^{-2Mr}.
  \end{equation}
  Hence, \begin{align*}
    \sum_{r=1}^{\infty}
    M^{K-1/4}(1+r)^K e^{-2Mr}
    &\le C_K M^{K-1/4}e^{-2M}
    \sum_{m=0}^{\infty}(2+m)^K e^{-2Mm}\\
    &\le C_K M^{K-1/4}e^{-2M}\longrightarrow0.
  \end{align*}
  Thus, \begin{equation}\label{eq:first-sum-o-rho-detail}
    \sum_{r=1}^{\lfloor\chi\sqrt{2M}\rfloor}
    C M^K(1+r)^K e^{-\mathcal J_M(r)}=o(\rho(M)).
  \end{equation}

  For $1\le r\le\chi\sqrt{2M}$ we have
  $(1+r)^K\le C M^{K/2}$, and there are $O(\sqrt M)$ terms. Therefore, \begin{align*}
    &\sum_{r=1}^{\lfloor\chi\sqrt{2M}\rfloor}
    C M^K(1+r)^K e^{-2\kappa M^{3/2}}
    \le C M^{K'}\sqrt M\,e^{-2\kappa M^{3/2}}.
  \end{align*}
  Dividing by $\rho(M)$ gives
  $
  C M^{K'+1/4}e^{-\kappa M^{3/2}}\longrightarrow0.
  $
  Consequently, the sum of the quadratic remainders is also $o(\rho(M))$.
  Combining this with~\eqref{eq:first-sum-o-rho-detail}, and replacing
  $r$ by $j-1$, yields
  \[
    \sum_{j=2}^{\lfloor\chi\sqrt{2M}\rfloor}
    \Pp(Y_0>M,Y_j>M)=o(\rho(M)).
  \]
  This completes the proof of Proposition~\ref{lem:separated}.
\end{proof}
\subsection{The trace formula and its four sectors}\label{SS:trace-sectors}
For the rest of this section, choose and fix $\chi\in((1/2)^{1/3},1)$. We consider integers $r$ with $1\le r\le\chi a$, where
\[
  a=\sqrt{2M},\qquad L=r+1,\qquad S=r+2,\qquad Q=Q_{M,1}.
\]
Each estimate below holds for all sufficiently large $M$, with the threshold depending only on $\chi$, and uniformly over this range of $r$. We also take $M$ large enough that $a\ge3/(1-\chi)$; then $a\ge S+1$. Constants denoted by $C_\chi,c_\chi>0$ may change from line to line. We write $\widetilde C_r=U_aC_rU_a^{-1}$ using the weighted-kernel convention of Section~\ref{sec:tail}.

Since $Q=P_1-\Lambda^M_{[0,1]}$,  we appeal to~\eqref{eq:Q-piecewise-L21} to write
\[
  Q(x,y)=Q^{\mathrm R}(x,y)+Q^{\mathrm F}(x,y),
\]
where
\[
  Q^{\mathrm R}(x,y)
  :=\mathbf 1_{\{x\le M,\ y\le M\}}p_1(2M-x-y)
\]
and
\[
  Q^{\mathrm F}(x,y)
  :=\mathbf 1_{\{x>M\ \mathrm{or}\ y>M\}}p_1(x-y).
\]
For $\alpha,\beta\in\{\mathrm R,\mathrm F\}$, put
\[
  C_r^{\alpha\beta}
  :=Q^{\alpha}P_rQ^{\beta}P_{-S}B_0.
\]
Then
\[
  C_r=C_r^{\mathrm{RR}}+C_r^{\mathrm{RF}}
  +C_r^{\mathrm{FR}}+C_r^{\mathrm{FF}},
  \qquad
  \Tr C_r=\sum_{\alpha,\beta}\Tr C_r^{\alpha\beta}.
\]

\begin{lemma}\label{L:interaction-trace-formula}
  Under these assumptions, $\widetilde C_r$ is trace class. For each $\alpha,\beta\in\{\mathrm R,\mathrm F\}$, the weighted trace of $C_r^{\alpha\beta}$ is the absolutely convergent integral
  \[
    \Tr C_r^{\alpha\beta}
    =\int_{\R^4}Q^\alpha(x,y)p_r(y-z)Q^\beta(z,w)
    (P_{-S}B_0)(w,x)\,dx\,dy\,dz\,dw.
  \]
  The sum of the four sector traces is $\Tr\widetilde C_r$.
\end{lemma}
\begin{proof}
  We first justify the trace integral. Put $L=r+1$ and
  \[
    A_r:=U_aQP_rQU_a\Phi,\qquad
    W_{S,a}:=\Phi^{-1}U_a^{-1}P_{-S}B_0U_a^{-1}.
  \]
  Since $0\le Q\le P_1$, we have $0\le QP_rQ\le QP_L$. Lemmas~\ref{lem:V} and~\ref{lem:UPhi} therefore imply that $A_r$ is Hilbert--Schmidt; Lemma~\ref{lem:W} gives the same property for $W_{S,a}$ because $a\ge S+1$. Thus, $\widetilde C_r=A_rW_{S,a}$ is trace class, and the trace formula for a product of Hilbert--Schmidt kernels gives
  \[
    \Tr\widetilde C_r=\iint_{\R^2} A_r(x,w)W_{S,a}(w,x)\,dx\,dw.
  \]
  The absolute value is integrable by Cauchy--Schwarz. Since the kernel of $A_r$ is nonnegative, Tonelli applied with $|W_{S,a}|$ also justifies expanding its two inner integrals; Fubini then gives the four-fold formula in the statement. The exponential and polynomial weights cancel. We write $\Tr C_r$ and $\Tr C_r^{\alpha\beta}$ for the corresponding weighted trace integrals. Each sector has the same justification by domination by $A_r$.
  This completes the proof of Lemma~\ref{L:interaction-trace-formula}.
\end{proof}
The letters $\mathrm R$ and $\mathrm F$ refer to the reflected and free parts of $Q$. We first examine the reflected--reflected sector, where the scalar minimization produces $\mathcal J_M(r)$. The sectors containing a free kernel are then compared with the same scalar function.

\subsection{The reflected--reflected sector}\label{SS:RR}
Put $\delta_{\max}=(a^2+S^2)/2$ and $\delta_*=2(r+\sqrt{a^2+r^2})$.
Define
\begin{align}\label{eq:Fdelta-mainproof}
  F(\delta):=\frac{\delta^2}{2}+ 2S(M-\delta)+\frac23S^3
  +\frac23(a^2-2\delta+S^2)^{3/2}, \quad \delta \in [0, \delta_{\mathrm max}].
\end{align}
\begin{lemma}\label{L:RR-trace}
  Under the standing assumptions, for all $0\le\delta\le\delta_{\max}$,
  \[
    F(\delta)\ge\mathcal J_M(r)+\tfrac12(\delta-\delta_*)^2.
  \]
  Moreover, $|\Tr C_r^{\mathrm{RR}}|\le C_\chi e^{-\mathcal J_M(r)}$.
\end{lemma}
\begin{proof}
  Recall that $(P_{-S}B_0)(w,x)=e^{-2S^3/3-S(w+x)}\Ai(w+x+S^2)$.
  We first estimate $\Tr C_r^{\mathrm{RR}}$.
  In the $\mathrm{RR}$ sector one has $x,y,z,w\le M$ and hence
  \[
    Q^{\mathrm R}(x,y)
    =(4\pi)^{-1/2}
    \exp\left\{-\frac{(x+y-2M)^2}{4}\right\},
  \]
  \[
    Q^{\mathrm R}(z,w)
    =(4\pi)^{-1/2}
    \exp\left\{-\frac{(z+w-2M)^2}{4}\right\}.
  \]
  We split this sector according to the sign of the Airy argument
  \[
    u:=w+x+S^2.
  \]
  We write
  \begin{align*}
    \Tr C_r^{\mathrm{RR}} =\Tr C_r^{\mathrm{RR},+}+ \Tr C_r^{\mathrm{RR},-},
  \end{align*}
  where
  \begin{align*}
    \Tr C_r^{\mathrm{RR},+}&
    =\int_{\substack{(x,y,z,w)\in \R^4\\w+x+S^2\geq0}}
    Q^{\mathrm R}(x,y)p_r(y-z)Q^{\mathrm R}(z,w)
    (P_{-S}B_0)(w,x)\,dx\,dy\,dz\,dw,\\
    \Tr C_r^{\mathrm{RR},-}&
    =\int_{\substack{(x,y,z,w)\in \R^4\\w+x+S^2<0}}
    Q^{\mathrm R}(x,y)p_r(y-z)Q^{\mathrm R}(z,w)
    (P_{-S}B_0)(w,x)\,dx\,dy\,dz\,dw.
  \end{align*}

  Using the Airy bound $|\Ai(u)|\le Ce^{-\frac23u^{3/2}}$, for $u\ge0$, we see that
  \begin{equation}\label{eq:Cr-trace-integral-mainproof}
    |\Tr C_r^{\mathrm{RR},+}|
    \le C r^{-1/2}
    \int_{\substack{x,y,z,w\le M\\ x+w+S^2\ge0}}
    \exp\{-\mathcal F(x,y,z,w)\}\,dx\,dy\,dz\,dw,
  \end{equation}
  where
  \begin{align*}
    \mathcal F(x,y,z,w)
    ={}&\frac{(x+y-2M)^2}{4}
    +\frac{(y-z)^2}{4r}
    +\frac{(z+w-2M)^2}{4} \notag\\
    &+S(x+w)+\frac23S^3
    +\frac23(x+w+S^2)^{3/2}.
    %\label{eq:F4-mainproof}
  \end{align*}
  Introduce symmetric and antisymmetric variables
  \begin{align}\label{eq:change}
    X=\frac{x+w}{2},\quad U=\frac{x-w}{2},\qquad
    Y=\frac{y+z}{2},\quad V=\frac{y-z}{2}.
  \end{align}
  The inverse transformation is
  \[
    x=X+U,\quad w=X-U,\qquad y=Y+V,\quad z=Y-V,
  \]
  and its Jacobian has absolute value $4$.  A direct calculation gives
  \[
    (x+y-2M)^2+(z+w-2M)^2
    =2(X+Y-2M)^2+2(U+V)^2
  \]
  and $(y-z)^2=4V^2$.  Hence, \begin{align*}
    \mathcal F
    =&\frac{(X+Y-2M)^2}{2}
    +\frac{(U+V)^2}{2}+\frac{V^2}{r}
    +2SX+\frac23S^3
    +\frac23(2X+S^2)^{3/2}.
    %\label{eq:Fsym-mainproof}
  \end{align*}
  Removing the half-space restrictions in $U,V$ can
  only increase the integral, and
  \begin{align}\label{eq:Gauss}
    \int_{\mathbb R^2}
    \exp\left\{-\frac{(U+V)^2}{2}-\frac{V^2}{r}\right\}\,dU\,dV
    =\sqrt{2\pi}\sqrt{\pi r}.
  \end{align}
  Thus, the factor $\sqrt r$ produced here cancels the $r^{-1/2}$ in~\eqref{eq:Cr-trace-integral-mainproof}. It follows that
  \begin{align*}
    |\Tr C_r^{\mathrm{RR},+}|
    &\le C \int_{-S^2/2}^M dX \int_{-\infty}^M dY
    \exp\left\{-\left[\frac{(X+Y-2M)^2}{2}+2SX\right.\right.\\
    &\hspace{47mm}\left.\left.{}+\frac23S^3
    +\frac23(2X+S^2)^{3/2}\right]\right\}.
  \end{align*}

  We make the change of variable $X=M-\delta$ with $\delta \in [0, \delta_{\mathrm max}]$.
  For fixed $X=M-\delta$, the
  $Y$ integral is
  \[
    \int_{-\infty}^{M}
    \exp\left\{-\frac{(X+Y-2M)^2}{2}\right\}\,dY
    =\int_{\delta}^{\infty}e^{-u^2/2}\,du.
  \]
  Recall the function $F$ in~\eqref{eq:Fdelta-mainproof}.

  Then we use the Gaussian integral bound $
  \int_{\delta}^{\infty}e^{-u^2/2}\,du
  \le C(1+\delta^{-1}\mathbf 1_{\{\delta\ge1\}})e^{-\delta^2/2}
  $ to derive that
  \begin{align}\label{eq:RR+}
    |\Tr C_r^{\mathrm{RR},+}|
    \le C
    \int_0^{\delta_{\max}}
    \bigl(1+\delta^{-1}\mathbf 1_{\{\delta\ge1\}}\bigr)
    e^{-F(\delta)}\,d\delta.
  \end{align}
  Put
  \[
    q(\delta)=\sqrt{a^2-2\delta+S^2}.
  \]
  Then
  \[
    F'(\delta)=\delta-2S-2q(\delta),
    \qquad
    F''(\delta)=1+\frac{2}{q(\delta)}>0.
  \]
  At the critical point
  $
  \delta=2(S+q),
  $
  substituting this into $q^2=a^2-2\delta+S^2$ yields
  \[
    (q+2)^2=a^2+(S-2)^2=a^2+r^2.
  \]
  Writing
  $
  R:=\sqrt{a^2+r^2},
  $
  we obtain
  $
  q_*=R-2$ and $\delta_*=2(r+R).
  $
  Because $r\le\chi a$ with $\chi<1$, we have $q_*>0$ and $\delta_*\in (0, \delta_{\mathrm max})$ for all sufficiently
  large $M$, uniformly in the stated range.  Strict convexity shows that $\delta_*$
  is the unique minimizer of $F$.  We substitute $\delta_*=2(r+R)$ into~\eqref{eq:Fdelta-mainproof} and use $S=r+2$ and some elementary calculations to deduce that
  \begin{align*}
    F(\delta_*)
    &=2(r+R)^2+Sa^2-4S(r+R)+\frac23S^3+\frac23(R-2)^3\\
    &=\frac23R^3+a^2r+\frac23r^3\\
    &=\frac23(a^2+r^2)^{3/2}+a^2r+\frac23r^3
    =\mathcal J_M(r).
  \end{align*}
  Thus, the exponent $\mathcal J_M(r)$ defined in~\eqref{eq:JMr} is obtained directly from the saddle
  point.

  Since $F''(\delta)=1+2/q(\delta)\ge1$ for $0\le\delta<\delta_{\max}$ and the unique
  minimizer $\delta_*=2(r+R)$ lies in the interior for all sufficiently large
  $M$, uniformly for $r\le\chi a$, Taylor's inequality and continuity at $\delta_{\max}$ give
  \begin{align}\label{eq:Fbound}
    F(\delta)\ge \mathcal J_M(r)+\frac12(\delta-\delta_*)^2, \qquad \text{for all $\delta \in [0, \delta_{\mathrm max}]$}.
  \end{align}
  Therefore, we conclude from~\eqref{eq:RR+} that
  \begin{align} \label{eq:CRR+}
    |\Tr C_r^{\mathrm{RR},+}|\le C e^{-\mathcal J_M(r)}  \int_\R e^{-\frac12(\delta-\delta_*)^2} d \delta= C' e^{-\mathcal J_M(r)}.
  \end{align}

  We estimate $\Tr C_r^{\mathrm{RR},-}$ next.  Using the Airy bound $|\Ai(u)|\leq C$ for $u<0$, we have
  \begin{equation*}%\label{eq:Cr-trace-integral-mainproof2}
    |\Tr C_r^{\mathrm{RR},-}|
    \le C r^{-1/2}
    \int_{\substack{x,y,z,w\le M\\ x+w+S^2<0}}
    \exp\{-\mathcal G(x,y,z,w)\}\,dx\,dy\,dz\,dw,
  \end{equation*}
  where
  \begin{align*}
    \mathcal G(x,y,z,w)
    ={}&\frac{(x+y-2M)^2}{4}
    +\frac{(y-z)^2}{4r}
    +\frac{(z+w-2M)^2}{4}
    +S(x+w)+\frac23S^3.
    %\label{eq:G4-mainproof}
  \end{align*}
  We perform the same change of variables as in~\eqref{eq:change}. Then remove the half-space restrictions in $U,V$ and appeal to~\eqref{eq:Gauss} to obtain
  \begin{equation*}%\label{eq:Cr-trace-integral-mainproof3}
    |\Tr C_r^{\mathrm{RR},-}|
    \le C \int_{-\infty}^{-S^2/2} d X \int_{-\infty}^M d Y \exp\left\{-\left[\frac{(X+Y-2M)^2}{2}
    +2SX+\frac23S^3
    \right]\right\}.
  \end{equation*}
  Recall $\delta_{\mathrm max}=\frac{a^2+S^2}{2} \geq 1$ and let $X=M-\delta$ again. It follows that
  \begin{align*}
    |\Tr C_r^{\mathrm{RR},-}|
    &\le C \int^{\infty}_{\delta_{\mathrm max}} d \delta \int^{\infty}_\delta d u \exp\left\{-\left[\frac{u^2}{2}
    +2S(M-\delta)+\frac23S^3\right]\right\}\\
    & \le C \int^{\infty}_{\delta_{\mathrm max}}  \exp\left\{-\left[\frac{\delta^2}{2}
    +2S(M-\delta)+\frac23S^3
    \right]\right\} d \delta,
  \end{align*}
  where the second inequality holds by the Gaussian integral bound $\int_{\delta}^\infty e^{-u^2/2}d u \leq e^{-\delta^2/2}$ for $\delta\geq 1$.
  Define
  \begin{equation*}%\label{eq:Gminus-mainproof}
    G(\delta)
    :=\frac{\delta^2}{2}+2S(M-\delta)+\frac23S^3,
    \qquad \delta\ge\delta_{\mathrm max}.
  \end{equation*}
  Its derivative is
  $
  G'(\delta)=\delta-2S.
  $
  Uniformly for $1\le r\le\chi a$, we have
  $\delta_{\mathrm max}=a^2/2+S^2/2>2S$ for all sufficiently
  large $M$.  Thus, $G$ is increasing on $[\delta_{\mathrm max}, \infty)$ and
  its minimum occurs at the boundary $\delta=\delta_{\mathrm max}$.  A direct substitution
  gives
  \begin{align*}
    G(\delta_{\mathrm max})
    &=\frac12\left(M+\frac{S^2}{2}\right)^2
    -\frac13S^3
    =\frac{a^4}{8}+\frac{a^2S^2}{4}+\frac{S^4}{8}-\frac{S^3}{3}.
    %\label{eq:Gminus-boundary-mainproof}
  \end{align*}
  Since
  $S=2+r\le2+\chi a$ uniformly for $1\le r\le\chi a$, there are constants $c_\chi>0$ and $a_0(\chi)$ such that
  \[
    G(\delta_{\mathrm max})\ge c_\chi a^4,
    \qquad a\ge a_0(\chi).
  \]
  Moreover,  since
  $
  G'(\delta)=\delta-2S.
  $
  we have
  \begin{align*}
    G'(\delta_{\mathrm max})
    &=M+\frac{S^2}{2} -2S
    =\frac{a^2}{2}+\frac{S^2}{2}-2S.
  \end{align*}
  It follows that there exist $c_\chi>0$ and
  $a_1(\chi)<\infty$ such that
  \[
    G'(\delta_{\mathrm max})\ge c_\chi a^2=c_\chi' M,
    \qquad a\ge a_1(\chi).
  \]
  Furthermore,  $G''(\delta)=1$, so $G'$ is increasing.  Consequently, for
  all $\delta\ge\delta_{\mathrm max}$,
  \begin{equation*}
    G(\delta)
    \ge G(\delta_{\mathrm max})+c_\chi' M(\delta-\delta_{\mathrm max}) \geq c_\chi M^2+c_\chi' M(\delta-\delta_{\mathrm max}).
    %\label{eq:Gminus-linear-growth}
  \end{equation*}
  Therefore, we obtain
  \begin{align}\label{eq:RR-}
    |\Tr C_r^{\mathrm{RR},-}| \leq C e^{-c_\chi M^2} \int_{\delta_{\mathrm max}}^\infty e^{-c'_\chi M(\delta-\delta_{\mathrm max})}d \delta \leq  C_\chi e^{-c_\chi M^2}.
  \end{align}

  Since $\mathcal J_M(r)=O_\chi(a^3)$, the negative-Airy contribution is absorbed into the stated bound for large $M$. This proves Lemma~\ref{L:RR-trace}.
  This completes the proof of Lemma~\ref{L:RR-trace}.
\end{proof}

\subsection{The mixed sectors}\label{SS:FR-RF}
\begin{lemma}\label{L:mixed-traces}
  Under the standing assumptions,
  \[
    \Tr C_r^{\mathrm{RF}}=\Tr C_r^{\mathrm{FR}},\qquad
    |\Tr C_r^{\mathrm{FR}}|\le C_\chi(1+r)^{1/2}e^{-\mathcal J_M(r)}.
  \]
\end{lemma}
\begin{proof}
  We next consider $C_r^{\mathrm{FR}}$.  Recall that
  \[
    C_r^{\mathrm{FR}}=Q^{\mathrm F}P_rQ^{\mathrm R}P_{-S}B_0,
    \qquad S=r+2,
  \]
  so that
  \begin{align}
    \Tr C_r^{\mathrm{FR}}
    ={}&\int_{\R^4}Q^{\mathrm F}(x,y)p_r(y-z)Q^{\mathrm R}(z,w)
    (P_{-S}B_0)(w,x)\,dx\,dy\,dz\,dw. \label{eq:FR-trace-start}
  \end{align}
  The free condition for the first $Q$ factor is the disjoint union
  \begin{align}\label{eq:D12}
    D_1:=\{x>M\},\qquad D_2:=\{x\le M,\ y>M\}.
  \end{align}
  We write the corresponding contributions as $I_1$ and $I_2$.
  \paragraph{The contribution $I_1$.}
  We first integrate over the unconstrained variable.
  Put
  $
  L:=r+1$ and hence $S=L+1$.
  Since $y$ is unrestricted on $D_1$, the semigroup property gives
  \[
    \int_{\R}p_1(x-y)p_r(y-z)\,dy=p_L(x-z).
  \]
  Consequently, with $\alpha=x-M$ and $\beta=M-w$, we have
  \begin{align*}%\label{eq:FR-I1-conv}
    I_1&=\int_{\substack{x>M\\z,w\le M}}
    p_L(x-z)p_1(2M-z-w)(P_{-S}B_0)(w,x)\,dx\,dz\,dw \nonumber\\
    &=\int_0^\infty d \alpha \int_{0}^\infty d \beta \, (P_{-S}B_0)(M-\beta,M+\alpha) \int_{-\infty}^M
    p_L(M+\alpha-z)p_1(M-z+\beta)\,dz \nonumber\\
    &= \int_0^\infty d \alpha \int_{0}^\infty d \beta \, (P_{-S}B_0)(M-\beta,M+\alpha) J(\alpha, \beta),
  \end{align*}
  where
  \[
    J(\alpha,\beta)
    =\int_{-\infty}^M
    p_L(M+\alpha-z)p_1(M-z+\beta)\,dz.
  \]

  By Lemma~\ref{lem:J},
  \begin{align*}
    J(\alpha,\beta)=p_{L+1}(\alpha-\beta)
    \Phi_{\rm G}\!\left(-\frac{\alpha+L\beta}{\sqrt{2L(L+1)}}\right),
  \end{align*}
  where $\Phi_{\rm G}$ is the standard normal distribution function.
  Using
  $\Phi_{\rm G}(-t)\le \frac12e^{-t^2/2}$ for $t\ge0$ gives
  \begin{align}\label{eq:boundJ}
    J(\alpha,\beta)
    &\le C S^{-1/2}
    \exp\left\{-\frac{(\alpha-\beta)^2}{4S}
    -\frac{(\alpha+L\beta)^2}{4LS}\right\}
    =C S^{-1/2}
    \exp\left\{-\frac{\alpha^2}{4L}-\frac{\beta^2}{4}\right\}.
  \end{align}
  On the positive-Airy part
  \[
    a^2+\alpha-\beta+S^2=x+w+S^2\ge0,
  \]
  we use the  Airy bound $|\Ai(u)|\leq C e^{-\frac23u^{3/2}}$ for $u\geq0$ and obtain
  \begin{equation}\label{eq:FR-I1-positive}
    |I_1^+|
    \le C S^{-1/2}
    \int_{\substack{\alpha,\beta\ge0\\a^2+\alpha-\beta+S^2\ge0}}
    e^{-H_1(\alpha,\beta)}\,d\alpha\,d\beta,
  \end{equation}
  where
  \begin{align*}
    H_1(\alpha,\beta)
    ={}&\frac{\alpha^2}{4L}+\frac{\beta^2}{4}
    +S(2M+\alpha-\beta)+\frac23S^3
    +\frac23(a^2+\alpha-\beta+S^2)^{3/2}. %\label{eq:H1-FR}
  \end{align*}

  Put $d=\alpha-\beta$ and recall the function $F$ defined in~\eqref{eq:Fdelta-mainproof}.  If $d\ge0$, then for fixed $d$ the quadratic
  term is minimized at $\beta=0$, $\alpha=d$, and the remaining part is increasing
  in $d$ because its derivative is
  $
  S+\sqrt{a^2+d+S^2}>0.
  $
  Thus, the minimum on $d\ge0$ occurs at $d=0$, and its value is at least
  $F(0)\ge\mathcal J_M(r)$.
  If $d<0$, write $d=-2\delta$ with $\delta>0$.  Under the constraints
  $\alpha,\beta\ge0$, the minimum of
  $\alpha^2/(4L)+\beta^2/4$ subject to $\alpha-\beta=-2\delta$ occurs at
  $\alpha=0$, $\beta=2\delta$, and equals $\delta^2$.  Hence, \begin{align*}
    \min_{\alpha-\beta=-2\delta}H_1(\alpha,\beta)
    ={}&\delta^2+2S(M-\delta)+\frac23S^3
    +\frac23(a^2-2\delta+S^2)^{3/2}\\
    ={}&F(\delta)+\frac{\delta^2}{2}
    \ge F(\delta)\ge\mathcal J_M(r).
  \end{align*}
  Therefore, \begin{equation}\label{eq:H1-lower}
    H_1(\alpha,\beta)\ge\mathcal J_M(r)
  \end{equation}
  throughout the positive-Airy part.

  We now estimate $I_1^+$. Denote its integration domain by
  \[
    \mathcal D_+
    :=\{(\alpha,\beta):\alpha\ge0,\ \beta\ge0,\
    a^2+\alpha-\beta+S^2\ge0\}.
  \]
  This is a closed convex set.  On its interior put
  $
  u=a^2+\alpha-\beta+S^2>0.
  $ Then the Hessian matrix of $H_1$ is given by
  \begin{align*}
    \nabla^2 H_1(\alpha,\beta)
    &=
    \begin{pmatrix}(2L)^{-1}&0\\0&1/2\end{pmatrix}
    +\frac1{2\sqrt u}
    \begin{pmatrix}1&-1\\-1&1\end{pmatrix}.
  \end{align*}
  Since $L=r+1\ge1$, it follows that for every vector $p\in \R^2$,
  \begin{align}\label{eq:Hessian}
    p^T \nabla^2 H_1(\alpha,\beta) p \geq \frac{1}{2L}|p|^2.
  \end{align}
  Let $p_0=(\alpha_0,\beta_0)$ be a constrained minimizer of $H_1$ on
  $\mathcal D_+$.  Existence follows from coercivity of the quadratic terms.
  Notice that one should not regard $p_0$ as an interior critical point:
  in fact, throughout the interior,
  \[
    \partial_\alpha H_1(\alpha,\beta)
    =\frac{\alpha}{2L}+S+\sqrt{a^2+\alpha-\beta+S^2}>0,
  \]
  so $\nabla H_1$ cannot vanish there.  Thus, the minimizer is necessarily a
  constrained boundary minimizer.

  First take an arbitrary interior point
  $p=(\alpha,\beta)\in\operatorname{int}(\mathcal D_+)$.  Since
  $\mathcal D_+$ is convex, the segment
  \[
    p_t:=p_0+t(p-p_0),\qquad 0<t\le1,
  \]
  lies in the interior for every $t>0$.  Define
  \[
    g(t):=H_1(p_t),\qquad 0\le t\le1.
  \]
  For $0<t\le1$, the function $g$ is twice differentiable and
  \[
    g''(t)
    =(p-p_0)^T\nabla^2H_1(p_t)(p-p_0)
    \ge \frac1{2L}|p-p_0|^2,
  \]
  where the inequality holds by~\eqref{eq:Hessian}. In addition,  we can compute that
  \begin{align*}
    \lim_{t\to0}\frac{g(t)-g(0)}{t}= \nabla H_1(p_0)\cdot(p-p_0).
  \end{align*}
  Since $g(t)\geq g(0)$ for all $0<t\leq 1$, it follows that $\nabla H_1(p_0)\cdot(p-p_0)\geq0$.
  For every $\varepsilon \in (0, 1)$
  Taylor's formula with integral remainder gives
  \begin{align}
    g(1)&= g(\varepsilon)+ g'(\varepsilon)(1-\varepsilon) +  \int_\varepsilon^1 (1-t)g''(t) d t \nonumber\\
    &\geq  g(\varepsilon)+ g'(\varepsilon)(1-\varepsilon) + \frac{|p-p_0|^2}{2L}\int_\varepsilon^1 (1-t) d t. \label{eq:ge}
  \end{align}
  From the expression of $g$, one can compute $g'(\varepsilon)$ and show that
  \[\lim_{\varepsilon\to0}g'(\varepsilon)= \nabla H_1(p_0)\cdot(p-p_0)\geq0.\]
  Hence, we let $\varepsilon \to 0$ in~\eqref{eq:ge} to obtain
  for every interior point $p$,
  \[
    H_1(p)
    \ge H_1(p_0)+\frac1{4L}|p-p_0|^2
    \ge \mathcal J_M(r)+\frac1{4L}|p-p_0|^2,
  \]
  where the last inequality uses~\eqref{eq:H1-lower}.

  Finally, if $p\in\partial\mathcal D_+$, choose interior points $p_n\to p$.
  Since $H_1$ is continuous on $\mathcal D_+$ and
  $|p_n-p_0|^2\to|p-p_0|^2$, passing to the limit yields
  \[
    H_1(\alpha,\beta)
    \ge \mathcal J_M(r)
    +\frac1{4L}
    \left| (\alpha,\beta)-(\alpha_0,\beta_0)\right|^2
  \]
  for every $(\alpha,\beta)\in\mathcal D_+$.
  Substituting this into~\eqref{eq:FR-I1-positive} and enlarging the domain of
  integration from $\mathcal D_+$ to $\mathbb R^2$ gives
  \begin{align*}
    |I_1^+|
    &\le C S^{-1/2}e^{-\mathcal J_M(r)}
    \int_{\mathbb R^2}
    \exp\left\{-\frac{|v|^2}{4L}\right\}\,dv
    = C S^{-1/2}(4\pi L)e^{-\mathcal J_M(r)}.
  \end{align*}
  Finally, $S=L+1$ and $L=r+1$, so
  $
  \frac{L}{\sqrt S}
  =\frac{r+1}{\sqrt{r+2}}
  \leq (1+r)^{1/2}.
  $
  It yields that
  \begin{equation}\label{eq:FR-I1-final}
    |I_1^+|\le C(1+r)^{1/2}e^{-\mathcal J_M(r)}.
  \end{equation}

  On the negative-Airy part one has
  \[
    a^2+\alpha-\beta+S^2<0,
    \qquad\text{hence}\qquad
    \beta>a^2+S^2+\alpha.
  \]
  Using $|\Ai(-u)|\le C$ for $u\ge0$, the backward Airy factor is therefore bounded by
  \[
    |(P_{-S}B_0)(M-\beta,M+\alpha)|
    \le C\exp\left\{-\frac23S^3-S(a^2+\alpha-\beta)\right\}.
  \]
  Combining this with the bound in~\eqref{eq:boundJ} gives
  \begin{equation}\label{eq:I1minus-double-int}
    |I_1^-|
    \le C S^{-1/2}
    \int_0^\infty\int_{a^2+S^2+\alpha}^\infty
    \exp\{-G_1(\alpha,\beta)\}\,d\beta\,d\alpha,
  \end{equation}
  where
  \begin{equation}\label{eq:G1minus}
    G_1(\alpha,\beta)
    :=\frac{\alpha^2}{4L}+\frac{\beta^2}{4}
    +S(a^2+\alpha-\beta)+\frac23S^3.
  \end{equation}
  For fixed $\alpha\ge0$, set
  $
  A(\alpha):=a^2+S^2+\alpha.
  $
  Then the domain of integration is $\beta\ge A(\alpha)$.  Since
  $
  \partial_\beta G_1(\alpha,\beta)=\frac\beta2-S
  $
  and
  \[
    A(\alpha)-2S
    =a^2+S^2+\alpha-2S
    =a^2+(S-1)^2-1+\alpha,
  \]
  we have, for all sufficiently large $a$,
  \[
    \partial_\beta G_1(\alpha,\beta)
    \ge \frac14a^2>0,
    \qquad \beta\ge A(\alpha),
  \]
  uniformly in $\alpha\ge0$.  Thus, for each fixed $\alpha$, the $G_1(\alpha,\beta)
  $ is
  increasing in $\beta$ throughout the allowed region and its minimum is attained
  at the boundary $\beta=A(\alpha)$.  Substituting this value into~\eqref{eq:G1minus} gives
  \begin{align}
    G_1(\alpha,A(\alpha))
    &=\frac{\alpha^2}{4L}
    +\frac{(a^2+S^2+\alpha)^2}{4}
    -\frac13S^3.\label{eq:G1minus-boundary}
  \end{align}
  In particular,
  \[
    G_1(\alpha,A(\alpha))
    \ge E_0,
    \qquad
    E_0:=\frac{(a^2+S^2)^2}{4}-\frac13S^3.
  \]
  Because $S=r+2\le \chi a+2$, we have $S^3\le C_\chi a^3$, while
  $(a^2+S^2)^2\ge a^4$.  Hence, there exist $c_\chi>0$ and $a_0(\chi)$ such that
  \begin{equation}\label{eq:E0-a4}
    E_0\ge c_\chi a^4,
    \qquad a\ge a_0(\chi).
  \end{equation}
  It remains to check that the integrations in $\beta$ and $\alpha$ introduce no
  exponential loss.  Write
  \[
    \beta=A(\alpha)+t,
    \qquad t\ge0.
  \]
  Since $G_1$ is quadratic in $\beta$,
  \begin{align*}
    G_1(\alpha,A(\alpha)+t)
    &=G_1(\alpha,A(\alpha))
    +\left(\frac{A(\alpha)}2-S\right)t+\frac{t^2}{4}.
  \end{align*}
  For all large $a$, we have uniformly in $\alpha\ge0$,
  $
  \frac{A(\alpha)}2-S\ge c a^2,
  $
  with a constant $c>0$.  Therefore, \begin{align*}
    \int_{A(\alpha)}^\infty e^{-G_1(\alpha,\beta)}\,d\beta
    &\le e^{-G_1(\alpha,A(\alpha))}
    \int_0^\infty e^{-c a^2t-t^2/4}\,dt
    \le \frac{C}{a^2}
    e^{-G_1(\alpha,A(\alpha))}.
  \end{align*}
  Next, expanding~\eqref{eq:G1minus-boundary},
  \[
    G_1(\alpha,A(\alpha))
    =E_0+\frac{a^2+S^2}{2}\alpha
    +\frac{\alpha^2}{4}+\frac{\alpha^2}{4L}
    \ge E_0+\frac{a^2}{2}\alpha.
  \]
  Hence, \[
  \int_0^\infty e^{-G_1(\alpha,A(\alpha))}\,d\alpha
  \le e^{-E_0}\int_0^\infty e^{-a^2\alpha/2}\,d\alpha
  =\frac{2}{a^2}e^{-E_0}.
\]
Substituting these estimates into~\eqref{eq:I1minus-double-int} yields
\[
  |I_1^-|
  \le C S^{-1/2}a^{-4}e^{-E_0}.
\]
Since $S\ge1$ and~\eqref{eq:E0-a4} holds uniformly in the Fredholm range
$r\le\chi a$, we conclude
\begin{equation}\label{eq:I1minus-final}
  |I_1^-|\le C_\chi e^{-c_\chi a^4}.
\end{equation}
Thus, the negative-Airy piece is superexponentially smaller on the
$M^{3/2}$ scale relevant for $e^{-\mathcal J_M(r)}$.

\smallskip
\paragraph{The contribution $I_2$.}
We now estimate $I_2$. Recall from~\eqref{eq:D12} that  $D_2=\{x\le M,\ y>M\}$.
Write
\[
  x=M-\alpha,\qquad z=M-\gamma,\qquad w=M-\beta,
  \qquad \alpha,\beta,\gamma\ge0.
\]
We first integrate over $y>M$:
\[
  K(\alpha,\gamma)
  :=\int_M^\infty p_1(x-y)p_r(y-z)\,dy.
\]
By the same calculation as in the proof of Lemma~\ref{lem:J},
we get
\begin{equation*}%\label{eq:FR-K-exact}
  K(\alpha,\gamma)
  =p_L(\gamma-\alpha)
  \Phi_{\rm G}\!\left(-\frac{r\alpha+\gamma}{\sqrt{2rL}}\right).
\end{equation*}
Therefore, using
$\Phi_{\rm G}(-t)\le \frac12e^{-t^2/2}$ for $t\ge0$ gives
\begin{align}\label{eq:Kbound}
  K(\alpha,\gamma)
  &\le C L^{-1/2}
  \exp\left\{-\frac{(\gamma-\alpha)^2}{4L}
  -\frac{(r\alpha+\gamma)^2}{4rL}\right\}
  =C L^{-1/2}
  \exp\left\{-\frac{\alpha^2}{4}-\frac{\gamma^2}{4r}\right\}.
\end{align}
The second $Q^{\mathrm R}$ factor  in~\eqref{eq:FR-trace-start} contributes
$
p_1(2M-z-w)=(4\pi)^{-1/2}e^{-(\gamma+\beta)^2/4}.
$
Thus, on the positive-Airy region $a^2-\alpha-\beta+S^2\ge0$,
\begin{equation}\label{eq:FR-I2-positive}
  |I_2^+|\le C L^{-1/2}
  \int_{\substack{\alpha,\beta,\gamma\ge0\\a^2-\alpha-\beta+S^2\ge0}}
  e^{-H_2(\alpha,\beta,\gamma)}\,d\alpha\,d\beta\,d\gamma,
\end{equation}
where
\begin{align*}
  H_2(\alpha,\beta,\gamma)
  ={}&\frac{\alpha^2}{4}+\frac{\gamma^2}{4r}
  +\frac{(\gamma+\beta)^2}{4}
  +S(2M-\alpha-\beta)+\frac23S^3 \notag\\
  &+\frac23(a^2-\alpha-\beta+S^2)^{3/2}. %\label{eq:H2-FR}
\end{align*}
We now estimate the integral in~\eqref{eq:FR-I2-positive} by integrating
first in $\gamma$.
Since $L=r+1$, completing the square gives the exact identity
\[
  \frac{\gamma^2}{4r}+\frac{(\gamma+\beta)^2}{4}
  =\frac{L}{4r}\left(\gamma+\frac{r\beta}{L}\right)^2
  +\frac{\beta^2}{4L}.
\]
Consequently,
\begin{align*}
  &\int_0^\infty
  \exp\left\{-\frac{\gamma^2}{4r}-\frac{(\gamma+\beta)^2}{4}\right\}
  \,d\gamma
  =
  2\sqrt{\frac rL}\,e^{-\beta^2/(4L)}
  \int_{\frac{\beta}{2}\sqrt{r/L}}^\infty e^{-u^2}\,du.
\end{align*}
Using the elementary Gaussian-tail estimate
$
\int_x^\infty e^{-u^2}\,du\le \frac{C}{1+x}e^{-x^2}$ for
$x\ge0
$
and the identity
$
\frac{\beta^2}{4L}+\frac{r\beta^2}{4L}=\frac{\beta^2}{4},
$
we obtain
\begin{align}\label{eq:gamma}
  \int_0^\infty
  \exp\left\{-\frac{\gamma^2}{4r}-\frac{(\gamma+\beta)^2}{4}\right\}
  \,d\gamma
  \le C\sqrt{\frac rL}\,e^{-\beta^2/4}.
\end{align}
Substituting this estimate into~\eqref{eq:FR-I2-positive} and using
$
L^{-1/2}\sqrt{\frac rL}=\frac{\sqrt r}{L}
=\frac{\sqrt r}{r+1}\le 1,
$
it gives
\[
  |I_2^+|
  \le C\int_{\substack{\alpha,\beta\ge0\\
  a^2-\alpha-\beta+S^2\ge0}}
  e^{-\widetilde H_2(\alpha,\beta)}\,d\alpha\,d\beta,
\]
where
\begin{align*}
  \widetilde H_2(\alpha,\beta)
  ={}&\frac{\alpha^2+\beta^2}{4}
  +S(2M-\alpha-\beta)+\frac23S^3
  +\frac23(a^2-\alpha-\beta+S^2)^{3/2}.
\end{align*}
Set
$
\delta:=\frac{\alpha+\beta}{2}$ and
$
v:=\frac{\alpha-\beta}{2}.
$
Then
\[
  \alpha=\delta+v,
  \qquad
  \beta=\delta-v,
  \qquad
  d\alpha\,d\beta=2\,d\delta\,dv,
\]
and
$
\frac{\alpha^2+\beta^2}{4}
=\frac{\delta^2+v^2}{2}.
$
Recalling the function $F$  defined in~\eqref{eq:Fdelta-mainproof}, we therefore have
\begin{equation}\label{eq:H2-exact-F-decomposition}
  \widetilde H_2(\alpha,\beta)
  =F(\delta)+\frac{v^2}{2}.
\end{equation}
The positive-Airy condition is precisely
$
0\le\delta\le\delta_{\max}:=\frac{a^2+S^2}{2}.
$
By the  estimate for $F$ in~\eqref{eq:Fbound}, it follows from~\eqref{eq:H2-exact-F-decomposition}  that
\[
  \widetilde H_2(\alpha,\beta)
  \ge \mathcal J_M(r)
  +\frac12(\delta-\delta_*)^2+\frac12v^2.
\]
We  now enlarge the transformed integration domain to all of $\mathbb R^2$.
Using the preceding bound, we obtain
\begin{align}
  |I_2^+|
  &\le C e^{-\mathcal J_M(r)}
  \int_{\mathbb R^2}
  \exp\left\{-\frac12(\delta-\delta_*)^2-\frac12v^2\right\}
  \,d\delta\,dv \nonumber\\
  &\le C e^{-\mathcal J_M(r)}. \label{eq:FR-I2-final}
\end{align}

We now treat the negative-Airy part of $D_2$.  Here
\[
  a^2-\alpha-\beta+S^2<0,
  \qquad\text{equivalently}\qquad
  \alpha+\beta>a^2+S^2.
\]
Since
$|\Ai(-u)|\le C$ for $u\ge0$ and
$
x+w=2M-\alpha-\beta=a^2-\alpha-\beta,
$
the backward Airy factor satisfies
\[
  |(P_{-S}B_0)(w,x)|
  \le C\exp\left\{-\frac23S^3-S(a^2-\alpha-\beta)\right\}.
\]
Combining this with the estimate for $K(\alpha,\gamma)$ in~\eqref{eq:Kbound} and the reflected
kernel $p_1(2M-z-w)$ obtained above gives
\begin{align*}
  |I_2^-|
  &\le C L^{-1/2}
  \int_{\substack{\alpha,\beta,\gamma\ge0\\
  \alpha+\beta>a^2+S^2}}
  \exp\{-G_2(\alpha,\beta,\gamma)\}
  \,d\alpha\,d\beta\,d\gamma, %\label{eq:I2minus-start}
\end{align*}
where
\[
  G_2(\alpha,\beta,\gamma)
  :=\frac{\alpha^2}{4}+\frac{\gamma^2}{4r}
  +\frac{(\gamma+\beta)^2}{4}
  +S(a^2-\alpha-\beta)+\frac23S^3.
\]
We integrate in $\gamma$ exactly as in the positive-Airy estimate.  The
bound in~\eqref{eq:gamma}
together with the fact $L^{-1/2}\sqrt{r/L}\le1$ yields
\begin{equation}\label{eq:I2minus-2d}
  |I_2^-|
  \le C\int_{\substack{\alpha,\beta\ge0\\
  \alpha+\beta>a^2+S^2}}
  \exp\{-\widetilde G_2(\alpha,\beta)\}
  \,d\alpha\,d\beta,
\end{equation}
with
\begin{align}\label{eq:G2}
  \widetilde G_2(\alpha,\beta)
  :=\frac{\alpha^2+\beta^2}{4}
  +S(a^2-\alpha-\beta)+\frac23S^3.
\end{align}
Introduce, as before,
$
\delta:=\frac{\alpha+\beta}{2}$ and
$v:=\frac{\alpha-\beta}{2}.
$
Then $d\alpha\,d\beta=2\,d\delta\,dv$ and
$
\frac{\alpha^2+\beta^2}{4}
=\frac{\delta^2+v^2}{2}.
$
The negative-Airy condition becomes
\[
  \delta>\delta_{\mathrm max},
  \qquad
  \delta_{\mathrm max}=\frac{a^2+S^2}{2}.
\]
Consequently, \[
\widetilde G_2(\alpha,\beta)
=G_{2,-}(\delta)+\frac{v^2}{2},
\]
where
\[
  G_{2,-}(\delta)
  :=\frac{\delta^2}{2}-2S\delta+Sa^2+\frac23S^3,
  \qquad \delta\ge\delta_{\mathrm max}.
\]
The function $G_{2,-}(\delta)$ is increasing throughout this region for large $a$.
Indeed,
\[
  G_{2,-}'(\delta)=\delta-2S
  \ge \delta_{\mathrm max}-2S
  =\frac{a^2+S^2-4S}{2}
  =\frac{a^2+r^2-4}{2},
\]
because $S=r+2$.  Hence, for all sufficiently large $a$,
\begin{equation}\label{eq:I2minus-slope}
  G_{2,-}'(\delta)\ge c a^2,
  \qquad \delta\ge\delta_{\mathrm max},
\end{equation}
with a constant $c>0$.  At the boundary,
\begin{align*}
  E_{2,0}
  :=G_{2,-}(\delta_{\mathrm max})
  &=\frac{(a^2+S^2)^2}{8}-\frac13S^3. %\label{eq:I2minus-boundary-cost}
\end{align*}
In the Fredholm range $r\le\chi a$ we have $S\le\chi a+2$, so
$S^3\le C_\chi a^3$, whereas $(a^2+S^2)^2\ge a^4$.  Therefore, for all
sufficiently large $a$ (depending at most on $\chi$),
\begin{equation}\label{eq:I2minus-a4}
  E_{2,0}\ge \frac18a^4-C_\chi a^3\ge c_\chi a^4.
\end{equation}
Moreover, by~\eqref{eq:I2minus-slope},
\[
  G_{2,-}(\delta)
  \ge E_{2,0}+c a^2(\delta-\delta_{\mathrm max}),
  \qquad \delta\ge\delta_{\mathrm max}.
\]
The constraints $\alpha,\beta\ge0$ imply $|v|\le\delta$; enlarging the
$v$-integration to all of $\mathbb R$ and using~\eqref{eq:I2minus-2d}, we get
\begin{align*}
  |I_2^-|
  &\le 2C e^{-E_{2,0}}
  \int_{\delta_{\mathrm max}}^\infty e^{-c a^2(\delta-\delta_{\mathrm \max})}\,d\delta
  \int_{\mathbb R}e^{-v^2/2}\,dv
  \le C a^{-2}e^{-E_{2,0}}.
\end{align*}
Combining this with~\eqref{eq:I2minus-a4} gives the desired estimate
\begin{equation}\label{eq:I2minus-final}
  |I_2^-|\le C_\chi e^{-c_\chi a^4}.
\end{equation}
Thus, this piece is again superexponentially small on the $a^3$ scale of
$\mathcal J_M(r)$.

Combining the four estimates in~\eqref{eq:FR-I1-final}, \eqref{eq:I1minus-final}, \eqref{eq:FR-I2-final}
\eqref{eq:I2minus-final},
and recalling that $a^4$ dominates
$\mathcal J_M(r)=O_\chi(a^3)$ in the Fredholm range, we conclude that
\begin{equation}\label{eq:FR-final-detailed}
  |\Tr C_r^{\mathrm{FR}}|
  \le C_\chi(1+r)^{1/2}e^{-\mathcal J_M(r)}.
\end{equation}
Under $(x,y,z,w)\mapsto(w,z,y,x)$ the two symmetric $Q$ kernels are interchanged, $p_r$ is even, and the Airy factor is unchanged. Thus, $\Tr C_r^{\mathrm{RF}}=\Tr C_r^{\mathrm{FR}}$, giving
\begin{equation}\label{eq:RF-final-detailed}
  |\Tr C_r^{\mathrm{RF}}|
  \le C_\chi(1+r)^{1/2}e^{-\mathcal J_M(r)}.
\end{equation}
This completes the proof of Lemma~\ref{L:mixed-traces}.
\end{proof}

\subsection{The free--free sector}\label{SS:FF}
\begin{lemma}\label{L:FF-trace}
  Under the standing assumptions,
  \[
    |\Tr C_r^{\mathrm{FF}}|\le C_\chi e^{-\mathcal J_M(r)}.
  \]
\end{lemma}
\begin{proof}
  It remains to estimate $\Tr C_r^{\mathrm{FF}}$.  Recall
  \[
    \Tr C_r^{\mathrm{FF}}
    =\int_{\mathbb R^4}
    Q^F(x,y)p_r(y-z)Q^F(z,w)(P_{-S}B_0)(w,x)
    \,dx\,dy\,dz\,dw,
  \]
  where
  \[
    Q^F(u,v)=\mathbf 1_{\{u>M\ \mathrm{or}\ v>M\}}p_1(u-v).
  \]
  It is convenient to decompose according to the two outer variables $x$ and
  $w$.  There are four pieces,
  \[
    \Tr C_r^{\mathrm{FF}}
    =T_{--}+T_{+-}+T_{-+}+T_{++},
  \]
  where the first sign records whether $x\le M$ or $x>M$, and the second sign
  records whether $w\le M$ or $w>M$.  We estimate these four pieces separately.

  \smallskip
  \paragraph{Both external variables below $M$.}
  We first estimate $T_{--}$.
  In this case the two free-kernel indicators force
  \[
    y>M,\qquad z>M.
  \]
  Write
  \[
    x=M-\alpha,\qquad w=M-\beta,
    \qquad y=M+\eta,\qquad z=M+\gamma,
  \]
  with $\alpha,\beta,\eta,\gamma\ge0$.  The Gaussian part is
  $
  p_1(\alpha+\eta)p_r(\eta-\gamma)p_1(\gamma+\beta).
  $
  Since
  \[
    (\alpha+\eta)^2\ge \alpha^2+\eta^2,
    \qquad
    (\gamma+\beta)^2\ge \gamma^2+\beta^2,
  \]
  we obtain
  \begin{align}\label{eq:GI}
    &\int_0^\infty\!\int_0^\infty
    p_1(\alpha+\eta)p_r(\eta-\gamma)p_1(\gamma+\beta)
    \,d\eta\,d\gamma \nonumber\\
    &\qquad\le
    C e^{-\alpha^2/4-\beta^2/4}
    \int_0^\infty e^{-\gamma^2/4}
    \left(\int_0^\infty p_r(\eta-\gamma)e^{-\eta^2/4}\,d\eta\right)d\gamma \nonumber\\
    &\qquad\le C e^{-\alpha^2/4-\beta^2/4}.
  \end{align}
  In the last step we used $e^{-\eta^2/4}\le1$ and
  $\int_{\mathbb R}p_r(\eta-\gamma)\,d\eta=1$, followed by the integrability of
  $e^{-\gamma^2/4}$.

  On the positive-Airy part,
  $
  a^2-\alpha-\beta+S^2\ge0,
  $
  the backward Airy factor therefore gives
  \begin{align*}
    |T_{--}^+|
    &\le C
    \int_{\substack{\alpha,\beta\ge0\\
    a^2-\alpha-\beta+S^2\ge0}}
    \exp\{-\widehat H_{--}(\alpha,\beta)\}
    \,d\alpha\,d\beta,
  \end{align*}
  where
  \[
    \widehat H_{--}(\alpha,\beta)
    :=\frac{\alpha^2+\beta^2}{4}
    +S(a^2-\alpha-\beta)+\frac23S^3
    +\frac23(a^2-\alpha-\beta+S^2)^{3/2}.
  \]
  Introduce again
  $
  \delta=\frac{\alpha+\beta}{2}$, $
  v=\frac{\alpha-\beta}{2}$.
  Then $d\alpha\,d\beta=2\,d\delta\,dv$ and, exactly as in~\eqref{eq:H2-exact-F-decomposition},
  \[
    \widehat H_{--}(\alpha,\beta)
    =F(\delta)+\frac{v^2}{2}.
  \]
  Using the estimate~\eqref{eq:Fbound} and
  enlarging the transformed domain to $\mathbb R^2$, we get
  \begin{equation}\label{eq:FF-Tmm-positive}
    |T_{--}^+|\le C e^{-\mathcal J_M(r)}.
  \end{equation}

  On the negative-Airy part, $\delta>\delta_{\mathrm max}:=(a^2+S^2)/2$.  Using only
  $|\Ai(-u)|\le C$ for $u\geq0$, the same Gaussian estimate in~\eqref{eq:GI} reduces  to
  \begin{equation*}%\label{eq:I2minus-2d2}
    |T_{--}^-|
    \le C\int_{\substack{\alpha,\beta\ge0\\
    \alpha+\beta>a^2+S^2}}
    \exp\{-\widetilde G_2(\alpha,\beta)\}
    \,d\alpha\,d\beta,
  \end{equation*}
  where $\widetilde G_2$ is given as in~\eqref{eq:G2}.
  Hence, the argument~\eqref{eq:I2minus-slope}--\eqref{eq:I2minus-a4} applies  and yields
  \begin{equation}\label{eq:FF-Tmm-negative}
    |T_{--}^-|\le C_\chi e^{-c_\chi a^4}.
  \end{equation}

  \smallskip
  \paragraph{One external variable above $M$.}
  We estimate $T_{+-}$, where $x>M$ and $w\le M$.
  Write
  \[
    x=M+\alpha,\qquad w=M-\beta,
    \qquad \alpha,\beta\ge0.
  \]
  Because $w\le M$, the second free factor forces $z>M$, while the first free
  factor is already active because $x>M$.  We may therefore integrate out $y$
  over the whole line:
  \[
    \int_{\mathbb R}p_1(x-y)p_r(y-z)\,dy=p_L(x-z),
    \qquad L=r+1.
  \]
  Writing $z=M+\gamma$, $\gamma\ge0$, the remaining Gaussian factor is
  $
  \int_0^\infty p_L(\alpha-\gamma)p_1(\gamma+\beta)\,d\gamma.
  $
  Since
  $
  p_1(\gamma+\beta)
  \le C e^{-\beta^2/4}e^{-\gamma^2/4}
  $
  and $e^{-\gamma^2/4}\le1$, we obtain the simple bound
  \begin{equation*}%\label{eq:FF-mixed-Gaussian}
    \int_0^\infty p_L(\alpha-\gamma)p_1(\gamma+\beta)\,d\gamma
    \le C e^{-\beta^2/4}.
  \end{equation*}
  Thus, on the positive-Airy region $a^2+\alpha-\beta+S^2\ge0$,
  \[
    |T_{+-}^+|
    \le C\int_{\substack{\alpha,\beta\ge0\\
    a^2+\alpha-\beta+S^2\ge0}}
    e^{-\Psi_+(\alpha,\beta)}\,d\alpha\,d\beta,
  \]
  where
  \[
    \Psi_+(\alpha,\beta)
    :=\frac{\beta^2}{4}
    +S(a^2+\alpha-\beta)+\frac23S^3
    +\frac23(a^2+\alpha-\beta+S^2)^{3/2}.
  \]
  Set $d=\alpha-\beta$.
  Recall the function $F$ defined in~\eqref{eq:Fdelta-mainproof} which satisfies
  \begin{align*}
    F(0)= Sa^2+\frac23S^3
    +\frac23(a^2+S^2)^{3/2} \geq \mathcal J_M(r).
  \end{align*}
  If $d\ge0$, then
  \[
    \Psi_+(\alpha,\beta)
    \ge F(0)+\frac{\beta^2}{4}+Sd
    \ge \mathcal J_M(r)+\frac{\beta^2}{4}+Sd.
  \]
  Changing variables from $(\alpha,\beta)$ to $(d,\beta)$ gives an integrable
  bound and hence a contribution at most $Ce^{-\mathcal J_M(r)}$.

  If $d<0$, write $d=-2\delta$ with $\delta>0$.  On the positive-Airy region, we have $0\leq \delta\leq \delta_{\mathrm max}:=\frac{a^2+S^2}{2}$. Since
  $\alpha=\beta-2\delta\ge0$, one has $\beta\ge2\delta$.  Put
  $\beta=2\delta+t$, $t\ge0$.  Then
  \begin{align*}
    \Psi_+(\alpha,\beta)
    &=\frac{(2\delta+t)^2}{4}
    +S(a^2-2\delta)+\frac23S^3
    +\frac23(a^2-2\delta+S^2)^{3/2}\\
    &=F(\delta)+\frac{\delta^2}{2}
    +\delta t+\frac{t^2}{4}.
  \end{align*}
  Therefore~\eqref{eq:Fbound} ensures
  \begin{align*}
    \Psi_+(\alpha,\beta)
    &\ge \mathcal J_M(r)
    +\frac12(\delta-\delta_*)^2
    +\frac{\delta^2}{2}+\delta t+\frac{t^2}{4}\\
    &\geq  \mathcal J_M(r)
    +\frac12(\delta-\delta_*)^2
    +\frac{t^2}{4}.
  \end{align*}
  The right-hand side is integrable in $(\delta,t)$, uniformly in $r$, and we
  conclude
  \begin{equation}\label{eq:FF-mixed-positive}
    |T_{+-}^+|\le C e^{-\mathcal J_M(r)}.
  \end{equation}

  On the negative-Airy part we use the bound $|\Ai(-u)|\le C$ for $u\geq0$ to get
  \[
    |T_{+-}^-|
    \le C\int_{\substack{\alpha,\beta\ge0\\
    a^2+\alpha-\beta+S^2<0}}
    e^{-\Psi_-(\alpha,\beta)}\,d\alpha\,d\beta,
  \]
  where
  \[
    \Psi_-(\alpha,\beta)
    :=\frac{\beta^2}{4}
    +S(a^2+\alpha-\beta)+\frac23S^3.
  \]
  Set $\alpha-\beta=-2\delta$. Then
  one necessarily has
  $\delta>\delta_{\mathrm max}=(a^2+S^2)/2$.  Again write $\beta=2\delta+t$, $t\ge0$.  It follows that
  \[
    \Psi_-(\alpha, \beta)=G_{+,-}(\delta)+\delta t+\frac{t^2}{4} \geq G_{+,-}(\delta)+ \frac{t^2}{4}
  \]
  with $
  G_{+,-}(\delta):=\delta^2-2S\delta+Sa^2+\frac23S^3.
  $
  Now
  \[
    G_{+,-}'(\delta)=2\delta-2S
    \ge a^2+S^2-2S\ge c a^2
  \]
  for all sufficiently large $a$, while at the boundary
  \[
    G_{+,-}(\delta_{\mathrm max})
    =\frac{(a^2+S^2)^2}{4}-\frac13S^3
    \ge c_\chi a^4.
  \]
  Thus, we have for all sufficiently large $a$
  \begin{align*}
    G_{+,-}(\delta) \geq c_\chi a^4 + c a^2(\delta-\delta_{\mathrm max}).
  \end{align*}
  Integrating first in $t$ and then in $\delta$ therefore gives
  \begin{equation}\label{eq:FF-mixed-negative}
    |T_{+-}^-|\le C_\chi e^{-c_\chi a^4}.
  \end{equation}
  By reversing the four Gaussian variables $(x,y,z,w)$, the same estimates hold
  for the other mixed piece:
  \begin{equation}\label{eq:FF-other-mixed}
    |T_{-+}^+|\le C e^{-\mathcal J_M(r)},
    \qquad
    |T_{-+}^-|\le C_\chi e^{-c_\chi a^4}.
  \end{equation}

  \smallskip
  \paragraph{Both external variables above $M$.}
  It remains to estimate $T_{++}$.
  Here both free indicators are automatically satisfied, so both intermediate
  variables may be integrated over the whole line.  By the semigroup property,
  \[
    \int_{\mathbb R^2}p_1(x-y)p_r(y-z)p_1(z-w)\,dy\,dz
    =p_S(x-w).
  \]
  Write $x=M+\alpha$, $w=M+\beta$, with $\alpha,\beta\ge0$.  The Airy argument is
  strictly positive, and therefore the Airy bound $|\Ai(u)| \leq C e^{-\frac23u^{3/2}}$ for $u\geq0$ yields that
  \begin{align*}
    |T_{++}|
    &\le C\int_0^\infty\!\int_0^\infty
    p_S(\alpha-\beta)
    \exp\{-R_+(\alpha+\beta)\}\,d\alpha\,d\beta,
  \end{align*}
  where
  \[
    R_+(t):=S(a^2+t)+\frac23S^3
    +\frac23(a^2+t+S^2)^{3/2}.
  \]
  With $t=\alpha+\beta$ and $v=\alpha-\beta$, the Jacobian is $1/2$ and
  $|v|\le t$.  Since $\int_{\mathbb R}p_S(v)\,dv=1$,
  \[
    |T_{++}|
    \le \frac C2\int_0^\infty e^{-R_+(t)}\,dt.
  \]
  Moreover, since
  $R_+'(t)=S+\sqrt{a^2+t+S^2}\ge S,
  $
  we have
  $
  R_+(t)\ge R_+(0)+St.
  $
  Finally, because $R_+(0)=F(0)\ge\mathcal J_M(r)$,  we obtain
  \begin{equation}\label{eq:FF-Tpp}
    |T_{++}|\le C S^{-1}e^{-F(0)}
    \le C e^{-\mathcal J_M(r)}.
  \end{equation}

  Combining~\eqref{eq:FF-Tmm-positive}, \eqref{eq:FF-Tmm-negative},
  \eqref{eq:FF-mixed-positive}, \eqref{eq:FF-mixed-negative},
  \eqref{eq:FF-other-mixed}, and~\eqref{eq:FF-Tpp}, and again absorbing the
  $e^{-c_\chi a^4}$ terms into the $e^{-\mathcal J_M(r)}$ scale, gives the following estimate
  \begin{equation}\label{eq:FF-final-detailed}
    |\Tr C_r^{\mathrm{FF}}|
    \le C_\chi e^{-\mathcal J_M(r)},
    \qquad 1\le r\le\chi a.
  \end{equation}
  This completes the proof of Lemma~\ref{L:FF-trace}.
\end{proof}

Combining~\eqref{eq:CRR+}, \eqref{eq:RR-}, \eqref{eq:FR-final-detailed}, \eqref{eq:RF-final-detailed} and~\eqref{eq:FF-final-detailed},  and  $\Tr\widetilde C_r=\Tr C_r$, gives the total trace bound
\begin{equation}\label{eq:traceC-main-direct}
  |\Tr\widetilde C_r|
  \le C_\chi(1+r)^{1/2}e^{-\mathcal J_M(r)},
  \qquad 1\le r\le\chi a.
\end{equation}

\subsection{Uniform trace-norm bounds}\label{SS:trace-norms}
\begin{lemma}\label{L:joint-trace-norms}
  Under the standing assumptions, there are $C,K>0$, depending only on $\chi$, such that
  \[
    \|\widetilde K_1\|_1+\|\widetilde K_2\|_1+
    \|\widetilde K_\cup\|_1+\|\widetilde C_r\|_1
    \le CM^K(1+r)^K e^{-2a^3/3}.
  \]
\end{lemma}
\begin{proof}
  We first estimate $\|\widetilde K_1\|_1$.  Recall~\eqref{eq:K12} and  Write
  \begin{align*}
    \widetilde K_1= (U_aQU_a\Phi) (\Phi^{-1}U_a^{-1}P_{-1}B_0U_a^{-1}).
  \end{align*}
  By Lemma~\ref{lem:V},
  \begin{equation}\label{eq:V-mainproof}
    \|U_aQU_a\Phi\|_2
    \le C M^K e^{a^2-a^3}
  \end{equation}
  for some $K>0$.
  Lemma~\ref{lem:W} implies that
  \begin{equation*}%\label{eq:W-mainproof}
    \|\Phi^{-1}U_a^{-1}P_{-1}B_0U_a^{-1}\|_2
    \le C e^{\frac13a^3-a^2}.
  \end{equation*}
  Thus, it follows that
  \begin{align}
    \|\widetilde K_1\|_1
    \le C M^K e^{-\frac23a^3}.
    \label{eq:K1-main-direct}
  \end{align}
  Moreover, recall~\eqref{eq:K12} and write
  \begin{align*}
    \widetilde K_2
    &=(U_aP_LU_a^{-1})
    (U_aQU_a\Phi)
    (\Phi^{-1}U_a^{-1}P_{-S}B_0U_a^{-1}).
  \end{align*}
  By Lemma~\ref{lem:2to2},
  \begin{equation}\label{eq:heat-main-direct}
    \|U_aP_LU_a^{-1}\|_{2\to2}\leq e^{La^2}.
  \end{equation}
  Using the ideal properties
  \[
    \|AT\|_1\le\|A\|_{2\to2}\|T\|_1,\qquad
    \|BC\|_1\le\|B\|_2\|C\|_2,
  \]
  we obtain
  \begin{align}
    \|\widetilde K_2\|_1
    &\le \|U_aP_LU_a^{-1}\|_{2\to2}
    \|U_aQU_a\Phi\|_2\|\Phi^{-1}U_a^{-1}P_{-S}B_0U_a^{-1}\|_2\notag\\
    &\le C M^K
    e^{La^2+a^2-a^3+\frac13a^3-Sa^2}
    =C M^K e^{-\frac23a^3},
    \label{eq:K2-main-direct}
  \end{align}
  where the second inequality holds by Lemmas~\ref{lem:V} and~\ref{lem:W}.

  It remains to estimate $\|\widetilde K_\cup\|_1$.  Recall~\eqref{eq:Kcup-def-main} (with $T=1$) and write
  \[
    K_\cup = D_rP_{-S}B_0, \qquad \text{with} \quad D_r:=P_1P_rP_1-\Lambda P_r\Lambda.
  \]
  Since $Q=P_1-\Lambda$ and all kernels are nonnegative,
  \[
    D_r=QP_rP_1+P_1P_rQ-QP_rQ,
    \qquad
    0\le D_r\le QP_L+P_LQ.
  \]
  Set
  \[
    \widetilde V_{D,r}:=U_aD_rU_a\Phi.
  \]
  For the term $P_LQ$,
  \[
    U_aP_LQU_a\Phi=(U_aP_LU_a^{-1})(U_aQU_a\Phi),
  \]
  so~\eqref{eq:heat-main-direct} and~\eqref{eq:V-mainproof} give
  \[
    \|U_aP_LQU_a\Phi\|_2
    \le C M^K e^{La^2+a^2-a^3}.
  \]
  For the term $QP_L$, by Lemma~\ref{lem:UPhi}
  \[
    \|U_aQP_LU_a\Phi\|_2
    \le C M^K(1+r)^K e^{La^2+a^2-a^3}.
  \]
  Consequently, \begin{equation}\label{eq:VD-main-direct}
    \|\widetilde V_{D,r}\|_2
    \le C M^K(1+r)^K e^{Sa^2-a^3}.
  \end{equation}
  Since
  \[
    U_aK_\cup U_a^{-1}=\widetilde V_{D,r} (\Phi^{-1}U_a^{-1}P_{-S}B_0U_a^{-1}),
  \]
  we obtain from~\eqref{eq:VD-main-direct} and Lemma~\ref{lem:W}
  \begin{equation}\label{eq:Kcup-main-direct}
    \|\widetilde K_\cup\|_1
    \le C M^K(1+r)^K e^{-\frac23a^3}.
  \end{equation}
  Finally, $
  \widetilde C_r=\widetilde K_1+\widetilde K_2-\widetilde K_\cup,
  $
  so~\eqref{eq:K1-main-direct}, \eqref{eq:K2-main-direct}, and~\eqref{eq:Kcup-main-direct} imply
  \begin{equation}\label{eq:crude-norm-main-direct}
    \|\widetilde K_1\|_1+\|\widetilde K_2\|_1
    +\|\widetilde K_\cup\|_1+\|\widetilde C_r\|_1
    \le C M^K(1+r)^K e^{-\frac23a^3}.
  \end{equation}
  This completes the proof of Lemma~\ref{L:joint-trace-norms}.
\end{proof}

\section{Mesoscopic anti-clustering}\label{sec:mes}
The preceding section controls joint block exceedances at separations in the range covered by Proposition~\ref{lem:separated} by a direct two-block Fredholm analysis. We now treat the remaining mesoscopic range. The argument is probabilistic: a quantitative modulus estimate reduces block exceedances to point exceedances on a fine grid, and Newman's inequality converts the cubic covariance decay into a threshold-sensitive joint-tail estimate. This avoids the need for translated trace-norm estimates outside the Fredholm range.

\subsection{Discretization of a block}\label{SS:discretization}
We first record a quantitative discretization estimate.  It is a direct consequence of the tightness estimate used in~\cite{BasuBhattacharjee2024} to prove convergence of the line-to-point LPP profile to Airy$_1$.

\begin{lemma}\label{lem:modulus}
  There exist constants $C,c,h_0>0$ such that, for $0<h\le h_0$,
  \begin{equation}\label{eq:modulus-tail}
    \Pp\left(
    \sup_{\substack{x,y\in[0,1], |x-y|\le 2h}}
    |\mathcal A_1(x)-\mathcal A_1(y)|>1
    \right)
    \le C h^{-1}\exp\left(-c h^{-1/2}\right).
  \end{equation}
\end{lemma}

\begin{proof}
  Lemma~5.1 of~\cite{BasuBhattacharjee2024}, combined with the compact-uniform weak convergence in their Theorem~1.4, implies the following local estimate: there are $C_0,c_0>0$ such that for $s\le1$ and $\ell$ sufficiently large,
  \begin{equation}\label{eq:BB-modulus-local}
    \Pp\left(
    \sup_{|t_1|,|t_2|\le s}
    |\mathcal A_1(t_1)-\mathcal A_1(t_2)|
    >C_0\ell s^{1/2}
    \right)
    \le C_0e^{-c_0\ell}.
  \end{equation}
  Indeed, in Lemma~5.1 of~\cite{BasuBhattacharjee2024} take their parameter $\varepsilon=2/3$; after the deterministic Airy$_1$ scaling, the threshold is a constant multiple of $\ell s^{1/2}$ and the right-hand side is $Ce^{-c\ell}$.  For each fixed $s$ and $\ell$, oscillation on the compact interval is a continuous functional, so its strict-exceedance event is open. The Portmanteau theorem therefore transfers the displayed upper bound to the Airy$_1$ limit.

  We now spell out the covering argument that turns~\eqref{eq:BB-modulus-local} into~\eqref{eq:modulus-tail}.  Partition $[0,1]$ into consecutive intervals of length $2h$,
  \[
    I_k=[2kh,2(k+1)h]\cap[0,1],
    \qquad 0\le k\le N-1,
  \]
  where, for $h$ sufficiently small,
  $
  N\le \frac{1}{2h}+1\le \frac{2}{h}.
  $
  By stationarity of the Airy$_1$ process, the local estimate~\eqref{eq:BB-modulus-local} applies, with the same constants, to every translate of an interval of radius $h$.  Take
  $
  s=h$ and
  $
  \ell=\frac{1}{2C_0\sqrt h}.
  $
  For $h$ small this value of $\ell$ is in the admissible large-$\ell$ range, and
  $
  C_0\ell h^{1/2}=\frac12.
  $
  Consequently, for every $k$,
  \[
    \Pp\left(
    \sup_{x,y\in I_k}|\mathcal A_1(x)-\mathcal A_1(y)|>\frac12
    \right)
    \le C_0\exp\left\{-\frac{c_0}{2C_0\sqrt h}\right\}.
  \]
  Let
  \[
    G_h:=\bigcap_{k=0}^{N-1}
    \left\{
    \sup_{x,y\in I_k}|\mathcal A_1(x)-\mathcal A_1(y)|\le\frac12
    \right\}.
  \]
  We claim that on $G_h$,
  \[
    \sup_{\substack{x,y\in[0,1], |x-y|\le2h}}
    |\mathcal A_1(x)-\mathcal A_1(y)|\le1.
  \]
  Indeed, suppose $x<y$ and $y-x\le2h$.  If $x$ and $y$ belong to the same interval $I_k$, then their increment is at most $1/2$.  Otherwise they belong to two consecutive intervals $I_k$ and $I_{k+1}$ (the only exceptional possibility is that they are boundary points at distance exactly $2h$, in which case both belong to the intervening closed interval).  Let
  $
  z=2(k+1)h
  $
  be the common endpoint of $I_k$ and $I_{k+1}$.  Since $x,z\in I_k$ and $z,y\in I_{k+1}$, the triangle inequality gives
  \[
    |\mathcal A_1(x)-\mathcal A_1(y)|
    \le |\mathcal A_1(x)-\mathcal A_1(z)|
    +|\mathcal A_1(z)-\mathcal A_1(y)|
    \le \frac12+\frac12=1.
  \]
  Hence, \[
  \left\{
  \sup_{\substack{x,y\in[0,1], |x-y|\le2h}}
  |\mathcal A_1(x)-\mathcal A_1(y)|>1
  \right\}
  \subset G_h^c.
\]
By the union bound,
\begin{align*}
  \Pp(G_h^c)
  &\le \sum_{k=0}^{N-1}
  \Pp\left(
  \sup_{x,y\in I_k}|\mathcal A_1(x)-\mathcal A_1(y)|>\frac12
  \right)\\
  &\le \frac{2C_0}{h}
  \exp\left\{-\frac{c_0}{2C_0}h^{-1/2}\right\}
  \le C h^{-1}e^{-c h^{-1/2}}.
\end{align*}
This is exactly~\eqref{eq:modulus-tail}, after renaming the constants.
This completes the proof of Lemma~\ref{lem:modulus}.
\end{proof}

\subsection{A threshold-sensitive covariance bound}\label{SS:newman}
We need the following high-threshold version of Newman's inequality.

\begin{lemma}\label{lem:tail-sensitive}
  Let $X=\mathcal A_1(0)$ and $Y=\mathcal A_1(d)$.  There are constants $C,u_0$ such that, for $u\ge u_0$ and $d\ge0$,
  \begin{equation}\label{eq:tail-sensitive}
    \Pp(X>u,Y>u)
    \le C\,\Pp(X>u)^2+C u\,\Cov(X,Y).
  \end{equation}
\end{lemma}

\begin{proof}
  Let $\varepsilon=u^{-1/2}$ and let $g=g_{u,\varepsilon}$ be the nondecreasing piecewise linear function satisfying
  \[
    g(x)=0\quad(x\le u-\varepsilon),\qquad
    g(x)=1\quad(x\ge u),\qquad
    \Lip(g)=\varepsilon^{-1}.
  \]
  Then $\mathbf 1_{(u,\infty)}\le g$. The Airy$_1$ process is associated (see~\cite{Pu2025}). Apply Newman's inequality~\cite[Lemma 3]{newman:80:normal} first to smooth mollifications of $g$, which converge uniformly to $g$, remain bounded between zero and one, and have derivative bounded by $\varepsilon^{-1}$. Bounded convergence then gives
  \[
    \Cov(g(X),g(Y))\le \varepsilon^{-2}\Cov(X,Y)=u\Cov(X,Y).
  \]
  We next derive the density estimate directly from the Painlev\'e--II representation of the GOE Tracy--Widom law in~\cite{BaikBuckinghamDiFranco2008}.  Let $F_1$ denote the standard GOE Tracy--Widom distribution function.  In the normalization used here,
  \[
    \Pp(X\le x)=F_1(2x),
    \qquad X=\mathcal A_1(0),
  \]
  so that, wherever the density exists,
  $
  f_X(x)=2F_1'(2x).
  $
  According to~\cite{BaikBuckinghamDiFranco2008},
  \[
    F_1(x)=F(x)E(x),
  \]
  where
  \[
    F(x)=\exp\left\{-\frac12\int_x^\infty R(s)\,ds\right\},
    \qquad
    E(x)=\exp\left\{-\frac12\int_x^\infty q(s)\,ds\right\},
  \]
  $q$ is the  solution of Painlev\'e~II,
  \[
    q''=2q^3+xq,
    \qquad
    q(x)\sim \Ai(x)\quad (x\to+\infty),
  \]
  and
  $
  R(x)=\int_x^\infty q(s)^2\,ds.
  $
  Differentiating the two exponential factors gives
  \[
    F'(x)=\frac12R(x)F(x),
    \qquad
    E'(x)=\frac12q(x)E(x),
  \]
  and therefore
  \begin{equation}\label{eq:F1-density-Painleve}
    F_1'(x)=\frac12\bigl(R(x)+q(x)\bigr)F_1(x).
  \end{equation}
  For large positive $x$, the  standard Airy asymptotics imply
  $
  0<q(x)\le Cx^{-1/4}e^{-\frac23x^{3/2}}.
  $
  Consequently, \begin{align*}
    R(x)
    &=\int_x^\infty q(s)^2\,ds \le C\int_x^\infty s^{-1/2}e^{-\frac43s^{3/2}}\,ds
    \le Cx^{-1}e^{-\frac43x^{3/2}},
  \end{align*}
  where the last inequality  is due to integration by parts formula.  Since $0<F_1(x)\le1$, \eqref{eq:F1-density-Painleve} yields, for all sufficiently large $x$,
  $
  F_1'(x)\le Cx^{-1/4}e^{-\frac23x^{3/2}}.
  $
  Returning to the Airy$_1$ normalization,
  \begin{align*}
    f_X(x)
    &=2F_1'(2x)
    \le Cx^{-1/4}
    \exp\left\{-\frac23(2x)^{3/2}\right\}
    =Cx^{-1/4}e^{-\kappa x^{3/2}}
  \end{align*}
  with $\kappa=\frac{4\sqrt2}{3}.$
  We now make this bound uniform over the short interval needed below.  If
  $v\in[u-\varepsilon,u]$ with $\varepsilon=u^{-1/2}$, then for all large $u$, we have
  $v\ge u/2$ and hence
  $
  v^{-1/4}\le 2^{1/4}u^{-1/4}.
  $
  Moreover, by the mean-value theorem,
  \[
    0\le u^{3/2}-v^{3/2}
    \le \frac32u^{1/2}(u-v)
    \le \frac32.
  \]
  Thus, \[
  e^{-\kappa v^{3/2}}
  =e^{-\kappa u^{3/2}}
  e^{\kappa(u^{3/2}-v^{3/2})}
  \le e^{3\kappa/2}e^{-\kappa u^{3/2}},
\]
and therefore
\begin{equation}\label{eq:fX-uniform-short-interval}
  \sup_{v\in[u-u^{-1/2},u]}f_X(v)
  \le C u^{-1/4}e^{-\kappa u^{3/2}}.
\end{equation}
The representation $F_1=FE$ above, with $q(x)\sim\Ai(x)$ and the smaller contribution from $\int_x^\infty R(s)\,ds$, gives $1-F_1(x)\sim\tfrac12\int_x^\infty\Ai(s)\,ds$; see also~\cite{BaikBuckinghamDiFranco2008}. Consequently,
$
\Pp(X>u)\asymp u^{-3/4}e^{-\kappa u^{3/2}}.
$
Since the interval $[u-\varepsilon,u]$ has length $\varepsilon=u^{-1/2}$, \eqref{eq:fX-uniform-short-interval} implies
\begin{align*}
  \Pp(u-\varepsilon<X\le u)
  &=\int_{u-\varepsilon}^{u}f_X(v)\,dv\\
  &\le C u^{-1/2}u^{-1/4}e^{-\kappa u^{3/2}}
  \le C\Pp(X>u).
\end{align*}
Consequently,
\[
  \Ee[g(X)]
  \le \Pp(X>u)+\Pp(u-\varepsilon<X\le u)
  \le C\Pp(X>u).
\]
Therefore, \[
\Pp(X>u,Y>u)
\le \Ee[g(X)g(Y)]
\le C\Pp(X>u)^2+u\Cov(X,Y),
\]
and the result follows after enlarging $C$.
This completes the proof of Lemma~\ref{lem:tail-sensitive}.
\end{proof}

\subsection{Summing the mesoscopic pairs}\label{SS:mesoscopic-sum}
\begin{proposition}\label{lem:mesoscopic}
  Fix $\chi\in((1/2)^{1/3},1)$.
  Let $M\to\infty$ and let $R_M=O(\sqrt M)$. Then
  \begin{equation}\label{eq:mesoscopic-sum}
    \sum_{j=\lceil\chi\sqrt{2M}\rceil}^{\lfloor R_M\rfloor}
    \Pp(Y_0>M,Y_j>M)
    =o(\rho(M))
  \end{equation}
  whenever the sum is nonempty.
\end{proposition}

\begin{proof}
  Choose $\vartheta\in(0,4/3)$ such that
  \begin{equation}\label{eq:theta-choice}
    \vartheta\chi^3>\frac23.
  \end{equation}
  By the sharp covariance asymptotic of~\cite{BasuBusaniFerrari2023},
  $\Cov(\mathcal A_1(0),\mathcal A_1(d))=\exp\{-(4/3+o(1))d^3\}$, we may choose $\vartheta<4/3$ satisfying~\eqref{eq:theta-choice} and then $d_0<\infty$ such that
  \begin{equation}\label{eq:cov-theta}
    \Cov(\mathcal A_1(0),\mathcal A_1(d))\le e^{-\vartheta d^3},
    \qquad d\ge d_0.
  \end{equation}

  Take the mesh size
  $
  h_M:=M^{-4},  n_{\rm mesh}:=\lceil h_M^{-1}\rceil+1=O(M^4)
  $
  and let $G_0$ and $G_j$ be the corresponding deterministic meshes of $[0,1]$ and $[j,j+1]$.  Define
  \[
    D_M(I):=\left\{
    \sup_{\substack{x,y\in I, |x-y|\le2h_M}}
    |\mathcal A_1(x)-\mathcal A_1(y)|>1
    \right\}.
  \]
  By Lemma~\ref{lem:modulus}, uniformly over translates $I$ of $[0,1]$,
  \begin{equation}\label{eq:DM-superexp}
    \Pp(D_M(I))\le C M^4e^{-cM^2}.
  \end{equation}
  If $Y_j>M$ and $D_M([j,j+1])$ does not occur, then some mesh point $x\in G_j$ satisfies $\mathcal A_1(x)>M-1$.  Hence, \begin{align}
    \Pp(Y_0>M,Y_j>M)
    &\le 2CM^4e^{-cM^2}
    +\sum_{x\in G_0}\sum_{y\in G_j}
    \Pp(\mathcal A_1(x)>M-1,\mathcal A_1(y)>M-1).\label{eq:block-to-grid}
  \end{align}

  By stationarity and Lemma~\ref{lem:tail-sensitive}, each summand is at most
  \[
    C p_*(M)^2+CM\,\Cov(\mathcal A_1(0),\mathcal A_1(|x-y|)),
    \qquad
    p_*(M):=\Pp(\mathcal A_1(0)>M-1).
  \]
  The GOE right-tail asymptotic gives
  \begin{equation}\label{eq:pstar}
    p_*(M)\le C M^{-3/4}
    \exp\left\{-\kappa M^{3/2}+C\sqrt M\right\}.
  \end{equation}
  For $j\ge\chi\sqrt{2M}$ and $x\in[0,1]$, $y\in[j,j+1]$, we have $|x-y|\ge j-1$, so~\eqref{eq:cov-theta} applies for large $M$.  Combining~\eqref{eq:block-to-grid}--\eqref{eq:pstar},
  \begin{align}
    \Pp(Y_0>M,Y_j>M)
    \le{}&CM^4e^{-cM^2}
    +CM^8p_*(M)^2
    +CM^9e^{-\vartheta(j-1)^3}.
    \label{eq:mesoscopic-point-bound}
  \end{align}

  There are $O(\sqrt M)$ indices in~\eqref{eq:mesoscopic-sum}.  The first contribution is superexponentially small relative to $\rho(M)$.  The second is bounded by a polynomial times
  $e^{-2\kappa M^{3/2}+C\sqrt M}=o(\rho(M))$.  For the third, write $a=\sqrt{2M}$.  Since
  $
  \kappa M^{3/2}=\frac23a^3
  $
  and $(j-1)^3\ge(\chi^3+o(1))a^3$, condition~\eqref{eq:theta-choice} gives
  \[
    M^{O(1)}\sqrt M\,e^{-\vartheta(j-1)^3}
    =o\left(e^{-2a^3/3}\right)=o(\rho(M)).
  \]
  This proves~\eqref{eq:mesoscopic-sum}.
  This completes the proof of Proposition~\ref{lem:mesoscopic}.
\end{proof}

\section{Proof of Theorem~\ref{thm:gumbel}}\label{sec:pf}

We prove Theorem~\ref{thm:gumbel} in this section. We first illustrate why we choose the centering constant $a_N$ as in~\eqref{eq:aN}.
\subsection{The extreme level}\label{SS:extreme-level}
Define
\begin{align}\label{eq:H}
  H(u):=\kappa u^{3/2}-\frac14\log u-\log C_*.
\end{align}
Proposition~\ref{lem:sharp-tail} gives
$
\Pp(Y_0>u)=e^{-H(u)}(1+o(1)).
$
The following lemma shows  why the coefficients in the centering constant $a_N$  are
$\gamma/9$ and
$(\gamma/9)\log(3/(2\pi^3))$.

\begin{lemma}\label{lem:aN}
  As $N\to\infty$,
  \[
    H(a_N)=\log N+o(1).
  \]
\end{lemma}
\begin{proof}
  Let $L_N:=\log N$.
  Write, for the moment,
  \[
    a_N=\gamma L_N^{2/3}+(A\log L_N+B)L_N^{-1/3}
    =\gamma L_N^{2/3}(1+\varepsilon_N),
  \]
  where
  \[
    \varepsilon_N
    :=\frac{A\log L_N+B}{\gamma L_N}
    =O\!\left(\frac{\log L_N}{L_N}\right).
  \]
  Since $\varepsilon_N\to0$, Taylor expansion gives
  \begin{align*}
    a_N^{3/2}
    &=\gamma^{3/2}L_N(1+\varepsilon_N)^{3/2}
    =\gamma^{3/2}L_N
    \left(1+\frac32\varepsilon_N+O(\varepsilon_N^2)\right).
  \end{align*}
  By the definition $\gamma=\kappa^{-2/3}$, we have
  $\kappa\gamma^{3/2}=1$.  Therefore, \begin{align*}
    \kappa a_N^{3/2}
    &=L_N+\frac{3}{2\gamma}(A\log L_N+B)
    +O(L_N\varepsilon_N^2)\\
    &=L_N+\frac{3A}{2\gamma}\log L_N
    +\frac{3B}{2\gamma}
    +O\!\left(\frac{(\log L_N)^2}{L_N}\right).
  \end{align*}
  Likewise,
  \begin{align*}
    \log a_N
    &=\log\gamma+\frac23\log L_N+\log(1+\varepsilon_N)\\
    &=\log\gamma+\frac23\log L_N
    +O\!\left(\frac{\log L_N}{L_N}\right).
  \end{align*}
  Substituting these two expansions into
  $
  H(u)=\kappa u^{3/2}-\frac14\log u-\log C_*
  $
  yields
  \begin{align}
    H(a_N)
    ={}&L_N
    +\left(\frac{3A}{2\gamma}-\frac16\right)\log L_N
    +\left(\frac{3B}{2\gamma}-\frac14\log\gamma-\log C_*\right)
    +o(1).
    \label{eq:HaN-direct-expansion}
  \end{align}
  Now take
  $
  A=\frac\gamma9$
  and $
  B=\frac\gamma9\log\left(\frac{3}{2\pi^3}\right).
  $
  Then
  $
  \frac{3A}{2\gamma}=\frac16,
  $
  so the coefficient of $\log L_N$ in~\eqref{eq:HaN-direct-expansion} vanishes.  It remains only to check the
  constant term.  Since
  $
  \gamma=\frac12\left(\frac32\right)^{2/3},
  $ $
  C_*=\frac{2^{1/4}}{\sqrt\pi},
  $
  we have
  \[
    \log\gamma=\frac23\log3-\frac53\log2,
    \qquad
    \log C_*=\frac14\log2-\frac12\log\pi.
  \]
  Consequently, \begin{align*}
    &\frac{3B}{2\gamma}-\frac14\log\gamma-\log C_*\\
    &\quad=
    \frac16\log\left(\frac{3}{2\pi^3}\right)
    -\frac14\left(\frac23\log3-\frac53\log2\right)
    -\left(\frac14\log2-\frac12\log\pi\right)
    =0.
  \end{align*}
  Thus, both correction terms cancel and
  \[
    H(a_N)=L_N+o(1)=\log N+o(1).
  \]
  This completes the proof of Lemma~\ref{lem:aN}.
\end{proof}

Recall the normalization constant $b_N$ defined in~\eqref{eq:bN}.
For fixed $z\in\R$, set
\begin{align}\label{eq:uN}
  u_N(z):=a_N+b_Nz.
\end{align}
\begin{lemma}\label{lem:uN}
  For every $z\in \R$,   we have $N\Pp\{Y_0>u_N(z)\}\to e^{-z}$ as $N\to\infty$.
\end{lemma}
\begin{proof}
  The function $H$ defined in~\eqref{eq:H} satisfies
  $
  H'(u)=2\sqrt2\,u^{1/2}-\frac1{4u},
  $
  so we have
  \[
    \frac1{H'(a_N)}
    \sim
    \frac{2\gamma}{3}(\log N)^{-1/3}=b_N.
  \]
  Hence, we have
  $
  H'(a_N)b_N\to 1
  $ as $N\to\infty$.
  Moreover,
  \[
    H''(u)=\sqrt2\,u^{-1/2}+\frac1{4u^2}.
  \]
  Because $a_N\sim\gamma(\log N)^{2/3}$ and
  $b_N=O((\log N)^{-1/3})$, for each fixed $z$ one has
  $u_N(z)\sim a_N$.  Taylor's theorem gives, for some
  $\xi_N$ between $a_N$ and $u_N(z)$,
  \begin{align*}
    H(u_N(z))
    &=H(a_N)+H'(a_N)b_Nz
    +\frac12H''(\xi_N)b_N^2z^2.
  \end{align*}
  Here
  $
  H'(a_N)b_Nz=z+o(1),
  $
  while
  \begin{align*}
    H''(\xi_N)b_N^2
    &=O\bigl((\log N)^{-1/3}\bigr)
    O\bigl((\log N)^{-2/3}\bigr)
    =O((\log N)^{-1})=o(1).
  \end{align*}
  Consequently, \[
  H(u_N(z))=\log N+z+o(1).
\]
Using $\Pp(Y_0>u)=e^{-H(u)}(1+o(1))$ at $u=u_N(z)$, we therefore obtain
\[
  N\Pp(Y_0>u_N(z))
  =\exp\{\log N-H(u_N(z))\}(1+o(1))
  \to e^{-z},
\]
as $N\to\infty$. This proves Lemma~\ref{lem:uN}.
This completes the proof of Lemma~\ref{lem:uN}.
\end{proof}

\subsection{Anti-clustering and blocking}\label{SS:blocking}
For the stationary sequence $(Y_j)$ and integers $m\ge1$, define
\[
  \alpha_Y(m):=\sup_{E\in\sigma(Y_j:j\le0),\ F\in\sigma(Y_j:j\ge m)}
  \left|\Pp(E\cap F)-\Pp(E)\Pp(F)\right|.
\]
For the Airy$_1$ process, $\alpha_{\mathcal A_1}(d)$ is defined using the fields generated by $\{\mathcal A_1(x):x\le0\}$ and $\{\mathcal A_1(x):x\ge d\}$. Since $Y_j$ depends on $[j,j+1]$, stationarity gives $\alpha_Y(m)\le\alpha_{\mathcal A_1}(m-1)$.

We are now ready to prove Theorem~\ref{thm:gumbel}.

\begin{proof}[{Proof of Theorem~\ref{thm:gumbel}}]
  Fix $z\in \R$ and let $\tau=e^{-z}$. Recall $u_N(z)$ defined in~\eqref{eq:uN} and put $p_N= \Pp(Y_0>u_N(z))$.
  We use the cubic-exponential strong-mixing estimate for the Airy$_1$ process proved in~\cite[Proposition 3.1]{ChenPu2026}:
  \[
    \alpha_{\mathcal A_1}(d)
    \le Ce^{-cd^3}.
  \]
  It follows that the sequence $\{Y_j\}_{j\in \mathbb{Z}}$ satisfies the following strong-mixing estimate
  \begin{equation}\label{eq:mixY}
    \alpha_Y(m)\le Ce^{-c(m-1)^3}.
  \end{equation}

  Choose
  \[
    q_N:=\left\lceil A(\log N)^{1/3}\right\rceil
  \]
  with $A$ sufficiently large that
  \begin{equation}\label{eq:mixsum}
    N\sum_{j\ge q_N}\alpha_Y(j)\longrightarrow0.
  \end{equation}
  Since
  $
  u_N(z)\asymp(\log N)^{2/3},
  $
  we have $q_N=O(\sqrt{u_N(z)})$. Fix
  $
  \chi\in\bigl((1/2)^{1/3},1\bigr).
  $ Put
  \[
    m_N:=\left\lfloor\chi\sqrt{2u_N(z)}\right\rfloor.
  \]
  We split the short-range sum at $m_N$. Lemma~\ref{lem:touching} gives
  \[
    N\Pp(Y_0>u_N,Y_1>u_N)\longrightarrow0.
  \]
  Proposition~\ref{lem:separated}, with $M=u_N(z)$, gives
  \[
    N\sum_{j=2}^{m_N}\Pp(Y_0>u_N,Y_j>u_N)\longrightarrow0,
  \]
  because $\Pp(Y_0>u_N)\sim\rho(u_N)$ by Proposition~\ref{lem:sharp-tail}, and $N\Pp(Y_0>u_N)\to\tau$ by Lemma~\ref{lem:uN}. Finally, since $q_N=O(\sqrt{u_N})$, Proposition~\ref{lem:mesoscopic} applies with $R_M=q_N$ and yields
  \[
    N\sum_{j=m_N+1}^{q_N}\Pp(Y_0>u_N,Y_j>u_N)\longrightarrow0.
  \]
  Therefore, \begin{equation}\label{eq:localDprime}
    N\sum_{j=1}^{q_N}\Pp(Y_0>u_N,Y_j>u_N)\longrightarrow0.
  \end{equation}

  For $j>q_N$, strong mixing implies
  \[
    \Pp(Y_0>u_N,Y_j>u_N)
    \le p_N^2+\alpha_Y(j).
  \]
  Thus, for each fixed integer $k\ge1$,
  \[
    N\sum_{j=q_N+1}^{\lfloor N/k\rfloor}
    \Pp(Y_0>u_N,Y_j>u_N)
    \le
    \frac{N^2}{k}p_N^2
    +N\sum_{j\ge q_N}\alpha_Y(j).
  \]
  Since $Np_N\to\tau$ by Lemma~\ref{lem:uN}, together with~\eqref{eq:mixsum}, we have
  \[
    \limsup_{N\to\infty}
    N\sum_{j=q_N+1}^{\lfloor N/k\rfloor}
    \Pp(Y_0>u_N,Y_j>u_N)
    \le\frac{\tau^2}{k}.
  \]
  Thus, by~\eqref{eq:localDprime}
  \begin{equation}\label{eq:Dprime}
    \lim_{k\to\infty}\limsup_{N\to\infty}
    N\sum_{j=1}^{\lfloor N/k\rfloor}
    \Pp(Y_0>u_N,Y_j>u_N)=0.
  \end{equation}

  We now perform a big-block/small-block argument with explicit indices. Fix an integer $k\ge1$ and put $b=\lfloor N/k\rfloor$. For $N$ large enough that $b>q_N$, take
  \[
    B_{m,N}:=\{(m-1)b,\ldots,mb-q_N-1\},\qquad 1\le m\le k.
  \]
  These consecutive cores have length $b-q_N=N/k+o(N)$ and consecutive cores are separated by $q_N$ discarded indices. The number of discarded indices, including the terminal remainder, is $kq_N+(N-kb)$. By Bonferroni's inequality,
  \[
    \sum_{j\in B_{m,N}}p_N
    -
    \sum_{\substack{i<j\\i,j\in B_{m,N}}}
    \Pp(Y_i>u_N,Y_j>u_N)
    \le
    \Pp\left\{\max_{j\in B_{m,N}}Y_j>u_N\right\}
    \le
    \sum_{j\in B_{m,N}}p_N.
  \]
  The first-order sum tends to $\tau/k$.  We now justify the order of the
  Bonferroni double sum.  Put
  $
  \ell_{m,N}:=|B_{m,N}|.
  $
  Since the $m$th long-block core has length $N/k+o(N)$, we may and do choose
  the blocking so that
  $
  \ell_{m,N}\le \left\lfloor\frac Nk\right\rfloor.
  $
  By stationarity, for $i<j$,
  \[
    \Pp(Y_i>u_N,Y_j>u_N)
    =
    \Pp(Y_0>u_N,Y_{j-i}>u_N).
  \]
  Grouping the pairs $(i,j)$ according to the separation $r=j-i$, there are
  exactly $\ell_{m,N}-r$ pairs in $B_{m,N}$ with separation $r$.  Hence, \begin{align*}
    \sum_{\substack{i<j\\ i,j\in B_{m,N}}}
    \Pp(Y_i>u_N,Y_j>u_N)
    &=
    \sum_{r=1}^{\ell_{m,N}-1}
    (\ell_{m,N}-r)
    \Pp(Y_0>u_N,Y_r>u_N)\\
    &\le
    \ell_{m,N}
    \sum_{r=1}^{\ell_{m,N}-1}
    \Pp(Y_0>u_N,Y_r>u_N)\\
    &\le
    \frac Nk
    \sum_{r=1}^{\lfloor N/k\rfloor}
    \Pp(Y_0>u_N,Y_r>u_N).
  \end{align*}
  Define
  \[
    D_{N,k}:=
    N\sum_{r=1}^{\lfloor N/k\rfloor}
    \Pp(Y_0>u_N,Y_r>u_N).
  \]
  Then the last display is bounded by $D_{N,k}/k$. The preceding local and mixing estimates give $0\le\limsup_{N\to\infty}D_{N,k}\le\tau^2/k$ for each fixed $k\ge1$. In particular, the following quantities are finite and, by~\eqref{eq:Dprime},
  \[
    \delta_k:=\limsup_{N\to\infty}D_{N,k}\longrightarrow0
    \qquad(k\to\infty).
  \]
  For each fixed sufficiently large $k$, Bonferroni gives, uniformly over the $k$ cores,
  \[
    1-\ell_{m,N}p_N\le \Pp(E_{m,N})
    \le 1-\ell_{m,N}p_N+D_{N,k}/k,
    \qquad E_{m,N}:=\{\max_{j\in B_{m,N}}Y_j\le u_N\}.
  \]
  Here $\ell_{m,N}p_N\to\tau/k$. Consequently, for every $\varepsilon>0$ and all sufficiently large $N$,
  \[
    1-\tau/k-\varepsilon\le\Pp(E_{m,N})
    \le1-\tau/k+\delta_k/k+\varepsilon.
  \]
  These bounds hold eventually in $N$ for each fixed $k$. Their remaining error is controlled as $k\to\infty$.

  Distinct cores are separated by $q_N$. By~\eqref{eq:mixY} and telescoping,
  \[
    \left|\Pp\Bigl(\bigcap_{m=1}^kE_{m,N}\Bigr)
    -\prod_{m=1}^k\Pp(E_{m,N})\right|
    \le(k-1)\alpha_Y(q_N)\longrightarrow0.
  \]
  If $G_N$ is the union of discarded indices, then $|G_N|=O(kq_N)$ and
  \[
    \Pp\{\exists j\in G_N:Y_j>u_N\}
    \le Ckq_Np_N\longrightarrow0
  \]
  for fixed $k$. For $k>\tau$, multiply the core bounds with $\varepsilon$ small enough that the lower base is positive. Let $N\to\infty$ through liminf and limsup, then let $\varepsilon\downarrow0$. We obtain
  \[
    (1-\tau/k)^k
    \le\liminf_{N\to\infty}\Pp\{\max_{0\le j<N}Y_j\le u_N\}
    \le\limsup_{N\to\infty}\Pp\{\max_{0\le j<N}Y_j\le u_N\}
    \le(1-\tau/k+\delta_k/k)^k.
  \]
  Both outer bounds tend to $e^{-\tau}$ as $k\to\infty$, since $\delta_k\to0$. This proves the limit $e^{-\tau}=e^{-e^{-z}}$.
  For real $N\to\infty$, put $n=\lfloor N\rfloor$. The explicit normalizations give
  \[
    \frac{a_N-a_n}{b_n}\to0,\qquad \frac{b_N}{b_n}\to1,
  \]
  and the same relations hold with $n+1$ in place of $n$.
  Since $\mathcal M_n\le\mathcal M_N\le\mathcal M_{n+1}$, the probability at
  threshold $u_N(z)$ is squeezed between the integer-interval probabilities at
  normalized levels $z\pm\varepsilon$, for every fixed $\varepsilon>0$ and all
  sufficiently large $N$. Let $N\to\infty$, then $\varepsilon\downarrow0$, and use
  continuity of the Gumbel distribution to obtain the asserted real-parameter limit.
  This completes the proof of Theorem~\ref{thm:gumbel}.
\end{proof}

\begin{remark}[The Gumbel limit from the estimates of Ferrari and Liu]\label{rem:FL}
  Ferrari and Liu~\cite{FL24} determined the fixed-level persistence exponent of the Airy$_1$ process. Their error bounds also give Theorem~\ref{thm:gumbel}. We sketch the deduction for comparison; it is not used in the proof above. Put
  \[
    P(c,L):=\Pp\{\mathcal A_1(x)\le c,\ 0\le x\le L\},\qquad
    \kappa_{\rm FL}(c):=-2\sum_{n\ge1}n^{-5/3}\Ai'(2n^{2/3}c).
  \]
  For $c\ge3/2$ and $L\ge1$, the proof of~\cite[Theorem 1.1]{FL24} in~\cite[Section 2.4]{FL24} gives
  \[
    \log P(c,L)=-L\kappa_{\rm FL}(c)-S(c)+E(c,L),
  \]
  where
  \[
    |S(c)|\le-6\log(1-2\beta),\qquad
    |E(c,L)|\le2e^{-4L^3/3}+\sum_{n\ge2}\frac{7^{n+1}\beta^n}{n}e^{-4L^3/(3n^2)},
  \]
  where the $n=1$ term is quoted from~\cite[Proposition 2.5]{FL24} and its intermediate estimates, the $n\ge2$ terms are bounded in~\cite[Section 2.4]{FL24}, and $\beta=\beta(c)$ is the quantity defined there with the choice $r=\sqrt{2c}$. The numerical constants in these bounds are absolute, and the dependence on $c$ enters only through $\beta(c)$. Recall that $\kappa=4\sqrt2/3$. The function $f(r)=r^3/3-2rc$ satisfies $f(\sqrt{2c})=-\kappa c^{3/2}$, and at $r=\sqrt{2c}$ we have $f(r-h)-f(r)=rh^2-h^3/3$. Taking $h=1/7$ gives
  \[
    \beta(c)\le2\exp\left\{-\kappa c^{3/2}+\frac{\sqrt{2c}}{49}\right\}.
  \]
  The simpler bound $\beta\le2e^{-\kappa c^{3/2}}$ in~\cite[Remark 2.4]{FL24} holds for moderate $c$: the second branch of their bound overtakes the first only for $c>578.4$. For large $c$, the second term in the definition of $\beta$ is the larger one, and the displayed bound covers all $c\ge3/2$. The right-hand side decreases in $c\ge1$, so $7\beta(c)\le\alpha_0:=14e^{-97\sqrt3/49}<0.46$ for all $c\ge3/2$. Beyond this uniform bound, only $\beta(c)\to0$ is used below. The bounds produced by the summation argument at the end of~\cite[Section 2.4]{FL24} increase in $\alpha=7\beta$ while the decay rate decreases, so $\alpha\le\alpha_0$ transfers them to constants independent of $c$; the argument is stated for large $L$, and the remaining range of $L$ is absorbed into $C$. Hence, $|E(c,L)|\le Ce^{-\delta L}$ with $C,\delta>0$ independent of $c\ge3/2$. Moreover, $S(c)\to0$ as $c\to\infty$.

  Now take $L=N$ and $c=u_N(z)$. The asymptotic $-2\Ai'(2c)=\rho(c)(1+O(c^{-3/2}))$ gives the term $n=1$, and the terms with $n\ge2$ contribute $O(c^{1/4}e^{-2\kappa c^{3/2}})$. The proof of Lemma~\ref{lem:uN} shows that $N\rho(u_N(z))=Ne^{-H(u_N(z))}\to e^{-z}$. Therefore, $N\kappa_{\rm FL}(u_N(z))\to e^{-z}$ and $P(u_N(z),N)\to e^{-e^{-z}}$.
\end{remark}

\appendix
\section{The continuum limit for two intervals}\label{S:continuum-appendix}

We now prove Lemma~\ref{lem:two-interval-determinant} by passing from finite meshes to the two constrained intervals.

\begin{proof}[Proof of Lemma~\ref{lem:two-interval-determinant}]
  Let $\bar P_M$ denote multiplication by $\mathbf 1_{(-\infty,M]}$.
  For $n\ge1$ choose a nested sequence of meshes of the two constrained
  intervals; for definiteness take dyadic $n=2^k$ and
  \[
    t_i=\frac{iT}{n},\qquad
    s_i=T+r+\frac{iT}{n},\qquad 0\le i\le n.
  \]
  Set
  \[
    E_n:=\bigcap_{i=0}^n\{\mathcal A_1(t_i)\le M\}
    \cap
    \bigcap_{i=0}^n\{\mathcal A_1(s_i)\le M\}.
  \]
  Since the meshes are nested and the Airy$_1$ process has continuous
  sample paths, $E_n$ decreases to
  \[
    E:=\{\mathcal A_1(x)\le M\text{ for all }
    x\in[0,T]\cup[T+r,2T+r]\}.
  \]
  Hence, we have
  \begin{equation}\label{eq:En-to-E}
    \Pp(E_n)\longrightarrow\Pp(E).
  \end{equation}

  We will apply Theorem~1 of~\cite{QuastelRemenik2012} to the ordered times
  \[
    t_0<t_1<\cdots<t_n<s_0<s_1<\cdots<s_n,
  \]
  with threshold $M$ at every sampled time. Let $\tau_0<\cdots<\tau_{2n+1}$
  denote the above ordered list. Then  equation (1.7) of~\cite{QuastelRemenik2012} contains the factor
  \[
    \bar P_M e^{-(\tau_0-\tau_1)\Delta}\bar P_M
    e^{-(\tau_1-\tau_2)\Delta}\cdots
    \bar P_M e^{-(\tau_{2n}-\tau_{2n+1})\Delta}\bar P_M
    e^{-(\tau_{2n+1}-\tau_0)\Delta}B_0.
  \]
  For our particular times this is
  \begin{align*}
    \mathcal T_n
    ={}&\bar P_M e^{-(t_0-t_1)\Delta}\bar P_M
    e^{-(t_1-t_2)\Delta}\cdots
    \bar P_M e^{-(t_{n-1}-t_n)\Delta}\bar P_M \notag\\
    &\quad\times e^{-(t_n-s_0)\Delta}\bar P_M
    e^{-(s_0-s_1)\Delta}\bar P_M\cdots
    e^{-(s_{n-1}-s_n)\Delta}\bar P_M
    \times e^{-(s_n-t_0)\Delta}B_0.
  \end{align*}
  Put $\delta=T/n$.  For
  $0\le i\le n-1$, we have
  $
  t_{i+1}-t_i=\delta$
  and $
  s_{i+1}-s_i=\delta.
  $
  Since $P_t=e^{t\Delta}$ and the times are increasing,
  \[
    e^{-(t_i-t_{i+1})\Delta}=e^{\delta\Delta}=P_\delta,
    \qquad
    e^{-(s_i-s_{i+1})\Delta}=P_\delta.
  \]
  Therefore, each constrained unit of length $T$ contributes the same
  finite-mesh operator
  \begin{equation}\label{eq:Lambda-n-twointerval}
    \Lambda_n
    :=\bar P_M P_\delta\bar P_M P_\delta\cdots
    P_\delta\bar P_M,
  \end{equation}
  which contains $n$ heat kernels and $n+1$ projections.

  Next, the last point of the first constrained interval and the first point
  of the second satisfy
  $
  s_0-t_n=(T+r)-T=r.
  $
  Hence, the factor crossing the gap is
  $
  e^{-(t_n-s_0)\Delta}=e^{r\Delta}=P_r.
  $
  There is no extra projection inside this factor: the event $E_n$ imposes
  no restriction at any time in the open gap $(T,T+r)$.  Finally,
  $
  s_n-t_0=(2T+r)-0=S,
  $
  so the closing factor in equation (1.7) of~\cite{QuastelRemenik2012} is
  $
  e^{-(s_n-t_0)\Delta}B_0=e^{-S\Delta}B_0=P_{-S}B_0.
  $
  Substituting these identities into the displayed formula for $\mathcal T_n$
  gives the exact factorization
  \[
    \mathcal T_n=\Lambda_nP_r\Lambda_nP_{-S}B_0.
  \]
  Thus, equation (1.7) of~\cite{QuastelRemenik2012}  yields
  \begin{equation}\label{eq:finite-twointerval-det}
    \Pp(E_n)
    =
    \det\left(I-B_0+\Lambda_nP_r\Lambda_nP_{-S}B_0\right).
  \end{equation}
  As in Remark~1.1 of~\cite{QuastelRemenik2012}, this Fredholm determinant is
  interpreted after the standard exponential conjugation.

  We next identify the fine-mesh limit.  We use the dyadic meshes chosen above (i.e. $n=2^k$), so that the mesh-point
  barrier events are nested.  We first derive an explicit formula for the kernel $\Lambda_n$  from its operator product in~\eqref{eq:Lambda-n-twointerval}.    Since
  \[
    (P_\delta f)(x)=\int_{\R}p_\delta(x-z)f(z)\,dz,
    \qquad
    (\bar P_Mf)(x)=\mathbf 1_{\{x\le M\}}f(x),
  \]
  the kernel of this product can be read off by applying the operators from
  right to left.  For instance, when $n=2$,
  \[
    \Lambda_2=\bar P_MP_\delta\bar P_MP_\delta\bar P_M,
  \]
  and therefore, for a test function $f$,
  \begin{align*}
    (\Lambda_2f)(x)
    &=\mathbf 1_{\{x\le M\}}
    \int_{\R}p_\delta(x-z)\mathbf 1_{\{z\le M\}}
    \left[\int_{\R}p_\delta(z-y)
    \mathbf 1_{\{y\le M\}}f(y)\,dy\right]dz\\
    &=\int_{\R}
    \left[
    \mathbf 1_{\{x\le M\}}\mathbf 1_{\{y\le M\}}
    \int_{-\infty}^{M}p_\delta(x-z)p_\delta(z-y)\,dz
    \right]f(y)\,dy.
  \end{align*}
  Thus, \[
  \Lambda_2(x,y)
  =\mathbf 1_{\{x\le M\}}\mathbf 1_{\{y\le M\}}
  \int_{-\infty}^{M}p_\delta(x-z)p_\delta(z-y)\,dz.
\]
The same computation for general $n$ gives, with $z_0=x$ and $z_n=y$,
\begin{align*}
  \Lambda_n(x,y)
  ={}&\mathbf 1_{\{x\le M\}}\mathbf 1_{\{y\le M\}}
  \int_{(-\infty,M]^{n-1}}
  \prod_{j=0}^{n-1}p_\delta(z_j-z_{j+1})
  \,dz_1\cdots dz_{n-1}.
\end{align*}

Let $B$ be Brownian motion with generator $\Delta$, so that
$
\mathbb P_x\{B_t\in dy\}=p_t(x-y)\,dy.
$
By the Markov property, the joint density of
$(B_{t_1},\ldots,B_{t_{n-1}},B_T)$, starting from $B_0=x$, is
\[
  \prod_{j=0}^{n-1}p_\delta(z_j-z_{j+1})
  \,dz_1\cdots dz_{n-1}\,dy.
\]
Consequently,
\begin{align*}
  &\mathbb P_x\left\{B_{t_j}\le M,
  j=1,\ldots,n-1,\ B_T\in dy\right\}\\
  &\qquad=
  \left[
  \int_{(-\infty,M]^{n-1}}
  \prod_{j=0}^{n-1}p_\delta(z_j-z_{j+1})
  \,dz_1\cdots dz_{n-1}
  \right]dy.
\end{align*}
On the other hand,
$
\mathbb P_x\{B_T\in dy\}=p_T(x-y)\,dy.
$
Dividing the preceding joint density by the endpoint density therefore gives
\begin{align*}
  &\frac{1}{p_T(x-y)}
  \int_{(-\infty,M]^{n-1}}
  \prod_{j=0}^{n-1}p_\delta(z_j-z_{j+1})
  \,dz_1\cdots dz_{n-1}\\
  &\qquad=
  \mathbb P^T_{x,y}
  \left\{b(t_j)\le M,\ j=1,\ldots,n-1\right\},
\end{align*}
where $\mathbb P^T_{x,y}$ denotes the Brownian-bridge law obtained by
conditioning $B_0=x$ on $B_T=y$, with the normalization corresponding to
$p_T$.    We have thus proved
\[
  \Lambda_n(x,y)
  =p_T(x-y)\,
  \mathbb P^{T}_{x,y}
  \left\{b(t_j)\le M\text{ for every }j=0,\ldots,n\right\}.
\]

Recall from~\cite[Theorem 4]{QuastelRemenik2012} that the kernel $\Lambda$ has the following representation
\[
  \Lambda(x,y)
  =p_T(x-y)\,
  \mathbb P^{T}_{x,y}
  \left\{b(t)\le M\text{ for every }0\le t\le T\right\}.
\]
Because the dyadic meshes are nested and Brownian-bridge paths are
continuous, it follows that for every
$x,y\in\R$,
\[
  \Lambda_n(x,y)\downarrow\Lambda(x,y).
\]
In particular,
\[
  0\le \Lambda_n-\Lambda\le P_T-\Lambda=:Q
\]
pointwise at the level of kernels. We write
\[
  (\Lambda_nP_r\Lambda_n)(x,w)
  =
  \iint_{\R^2}
  \Lambda_n(x,y)p_r(y-z)\Lambda_n(z,w)\,dy\,dz,
\]
and
\[
  (\Lambda P_r\Lambda)(x,w)
  =
  \iint_{\R^2}
  \Lambda(x,y)p_r(y-z)\Lambda(z,w)\,dy\,dz.
\]
Since $0\le\Lambda_n,\Lambda\le P_T$ and
\[
  \iint_{\R^2}
  p_T(x-y)p_r(y-z)p_T(z-w)\,dy\,dz
  =p_{2T+r}(x-w)
\]
by the semigroup property of heat kernel,
the dominated convergence theorem therefore yields that for every $x,w$,
\[
  (\Lambda_nP_r\Lambda_n)(x,w)
  \longrightarrow
  (\Lambda P_r\Lambda)(x,w).
\]
We now upgrade the pointwise kernel convergence to the
trace-norm convergence.
Set
\[
  D_n:=P_S-\Lambda_nP_r\Lambda_n,
  \qquad
  D:=P_S-\Lambda P_r\Lambda.
\]
Then
\[
  D_n-D
  =-(\Lambda_nP_r\Lambda_n-\Lambda P_r\Lambda),
\]
and the difference in parentheses has the exact telescoping decomposition
\[
  \Lambda_nP_r\Lambda_n-\Lambda P_r\Lambda
  =(\Lambda_n-\Lambda)P_r\Lambda_n
  +\Lambda P_r(\Lambda_n-\Lambda).
\]
All kernels in this display are nonnegative.  Using
$0\le\Lambda_n-\Lambda\le Q$ and
$0\le\Lambda_n,\Lambda\le P_T$, we obtain the pointwise domination
\begin{align*}
  |D_n-D|
  &\le QP_rP_T+P_TP_rQ
  =QP_{r+T}+P_{T+r}Q
  =:D_r^{\rm dom}.
\end{align*}
Moreover, $D=QP_{r+T}+\Lambda P_rQ$ and $\Lambda_n\ge\Lambda$ give
\[
  0\le D_n\le D\le D_r^{\rm dom}.
\]

Put $a=\sqrt{2M}$ and recall $U_a$ and $\Phi$ defined in~\eqref{eq:UPhi}.
Set
\[
  V_n:=U_aD_nU_a\Phi,
  \qquad
  V:=U_aDU_a\Phi,
\]
and
\[
  W_{S,a}:=\Phi^{-1}U_a^{-1}P_{-S}B_0U_a^{-1}.
\]
It follows that
\[
  U_aD_nP_{-S}B_0U_a^{-1}=V_nW_{S,a},
  \qquad
  U_aDP_{-S}B_0U_a^{-1}=VW_{S,a}.
\]
Lemma~\ref{lem:W} proves
$
\|W_{S,a}\|_2<\infty.
$
Moreover, we write
\begin{align*}
  U_aQP_{r+T}U_a \Phi= U_aQU_a\Phi \Phi^{-1}U_a^{-1} P_{r+T}U_a \Phi.
\end{align*}
The estimate~\eqref{eq:B-Schur-final}, proved in Lemma~\ref{lem:UPhi}, gives
\[
  \|\Phi^{-1}U_a^{-1}P_{r+T}U_a\Phi\|_{2\to2}<\infty.
\]
Combining this with Lemma~\ref{lem:V}, we obtain
$\|U_aQP_{r+T}U_a\Phi\|_2<\infty$. Similarly,
\[
  U_aP_{T+r}QU_a\Phi
  =U_aP_{T+r}U_a^{-1}U_aQU_a\Phi.
\]
Lemma~\ref{lem:2to2} gives
\[
  \|U_aP_{T+r}U_a^{-1}\|_{2\to2}<\infty.
\]
Together with Lemma~\ref{lem:V}, this yields
\[
  \|U_aP_{T+r}QU_a\Phi\|_2<\infty.
\]
Consequently, we obtain
\[
  \|U_aD_r^{\rm dom}U_a\Phi\|_2<\infty.
\]
The same domination proves that $V_n$ and $V$ are themselves Hilbert--Schmidt, uniformly in $n$.
On the other hand, the pointwise convergence just proved implies
$D_n(x,w)\to D(x,w)$, while the preceding domination gives
\[
  \bigl|(U_a(D_n-D)U_a\Phi)(x,w)\bigr|
  \le
  \bigl|(U_aD_r^{\rm dom}U_a\Phi)(x,w)\bigr|.
\]
The square of the right-hand side is integrable on $\R^2$.  Dominated
convergence in $L^2(\R^2)$ therefore yields
\[
  \|V_n-V\|_2
  =
  \|U_a(D_n-D)U_a\Phi\|_2
  \longrightarrow0.
\]
Using the ideal property $\|AB\|_1\le\|A\|_2\|B\|_2$, we conclude that
\begin{align*}
  &\left\|
  U_a\Bigl[
  \Lambda_nP_r\Lambda_nP_{-S}B_0
  -
  \Lambda P_r\Lambda P_{-S}B_0
  \Bigr]U_a^{-1}
  \right\|_1\\
  &\qquad\qquad=
  \|(V_n-V)W_{S,a}\|_1
  \le
  \|V_n-V\|_2\,\|W_{S,a}\|_2
  \longrightarrow0.
\end{align*}
Define the full trace-class perturbations
\[
  \widetilde K_n:=V_nW_{S,a}
  =U_a(B_0-\Lambda_nP_r\Lambda_nP_{-S}B_0)U_a^{-1},
  \qquad \widetilde K:=VW_{S,a}.
\]
The convergence just proved is precisely
\begin{equation}\label{eq:twointerval-operator-limit}
  \|\widetilde K_n-\widetilde K\|_1\longrightarrow0.
\end{equation}
The trace-class bounds apply to the complete weighted differences containing $B_0$.

Finally, using the inequality
\[
  |\det(I+A)-\det(I+B)|
  \le
  \|A-B\|_1
  \exp\{1+\|A\|_1+\|B\|_1\},
\]
with $A=-\widetilde K_n$ and $B=-\widetilde K$, we obtain~\eqref{eq:two-interval-continuum} by combining~\eqref{eq:finite-twointerval-det},
\eqref{eq:twointerval-operator-limit}, and~\eqref{eq:En-to-E}.
This completes the proof of Lemma~\ref{lem:two-interval-determinant}.
\end{proof}

\section{Auxiliary operator estimates}\label{S:operator-appendix}
\subsection{Airy integrals and Hilbert--Schmidt norms}
\begin{lemma}\label{lem:airy-laplace}
  There is a constant $C>0$ such that for every $\beta\ge1$,
  \begin{equation}\label{eq:airy-laplace}
    \int_{\R}e^{2\beta u}\Ai(u)^2\,du
    \le C e^{\frac23\beta^3}.
  \end{equation}
\end{lemma}

\begin{proof}
  For $u\le0$ we use the fact $|\Ai(u)|\le C$  to give a contribution bounded by a constant.  For $u\ge0$, put $u=v^2$ and use
  $
  |\Ai(v^2)|^2\le C(1+v^2)^{-1/2}e^{-\frac43v^3}.
  $
  Since $2v(1+v^2)^{-1/2}\le2$, the positive part is bounded by
  \begin{align*}
    &C\int_0^\infty
    \exp\left\{2\beta v^2-\frac43v^3\right\}dv\\
    &\qquad\qquad= Ce^{\frac23\beta^{3}} \int_0^\infty e^{-\frac23(v-\beta)^2(\beta+2v)}d v
    \leq C e^{\frac23\beta^{3}}\int_0^\infty e^{-\frac23(v-\beta)^2}d v\\
    &\qquad\qquad \leq C e^{\frac23\beta^{3}}\int_\R e^{-\frac23v^2}d v
    =C'e^{\frac23\beta^{3}}.
  \end{align*}
  This proves~\eqref{eq:airy-laplace}.
  This completes the proof of Lemma~\ref{lem:airy-laplace}.
\end{proof}

Fix $T>0$ and $\lambda\ge T+1$. Let $V_{T, \lambda}$ and $W_{T, \lambda}$ be two kernels defined as in~\eqref{eq:V-L21}
and~\eqref{eq:W-L21} respectively, with $a$ replaced by $\lambda$. The following two lemmas give bounds on their Hilbert--Schmidt norms.

\begin{lemma} \label{lem:W}
  There is an absolute constant $C>0$ such that, for every $T>0$ and every $\lambda\ge T+1$,
  \begin{equation}\label{eq:W-HS}
    \|W_{T,\lambda}\|_2^2
    \le C
    e^{\frac23\lambda^3-2T\lambda^2}.
  \end{equation}
\end{lemma}
\begin{proof}
  Using
  $
  (P_{-T}B_0)(z,y)
  =e^{-2T^3/3-T(z+y)}\Ai(z+y+T^2)
  $, we write
  \begin{align}
    \|W_{T,\lambda}\|_2^2&=
    \int_{\R^2} \phi^{-2}(z) e^{2\lambda(z+y)}e^{-4T^3/3-2T(z+y)}\Ai(z+y+T^2)^2dy dz \nonumber\\
    &=\|\phi^{-1}\|_2^2
    e^{\frac23T^3-2\lambda T^2}
    \int_{\R}e^{2(\lambda-T) u}\Ai(u)^2\,du\notag\\
    &\le C
    \exp\left\{-\frac43T^3-2(\lambda-T) T^2+\frac23(\lambda-T)^3\right\}, \nonumber
  \end{align}
  where the inequality is due to Lemma~\ref{lem:airy-laplace}.
  A direct expansion yields
  \begin{equation*}
    -\frac43T^3-2(\lambda-T)T^2+\frac23(\lambda-T)^3
    =\frac23\lambda^3-2T\lambda^2,
  \end{equation*}
  which implies~\eqref{eq:W-HS}.
  This completes the proof of Lemma~\ref{lem:W}.
\end{proof}

\begin{lemma} \label{lem:V}
  Fix $T>0$. There exists a constant $C_T>0$ such that, for every $M>0$ and every $\lambda\ge T+1$, writing $a=\sqrt{2M}$, we have
  \begin{equation*}\label{eq:W-HS2}
    \|V_{T,\lambda}\|_2^2
    \le C_T(1+M)^2(1+\lambda^3)
    e^{2T\lambda^2-2a^2\lambda}.
  \end{equation*}
\end{lemma}
\begin{proof}
  In the reflected region $x,z\le M$, we see from~\eqref{eq:Q-piecewise-L21} that
  \[
    Q_{M,T}(x,z)
    =\frac1{\sqrt{4\pi T}}
    e^{-(x-z)^2/(4T)}e^{-(M-x)(M-z)/T}.
  \]
  Hence,
  \begin{align*}
    \iint_{x,z\leq M}|V_{T,\lambda}(x,z)|^2 dx dz=\frac{1}{4\pi T}\iint_{x,z\leq M} (1+|z|^2)e^{-2\lambda(x+z)}
    e^{-\frac{(x-z)^2}{2T}}e^{-\frac{2(M-x)(M-z)}{T}}dx dz.
  \end{align*}
  Using change of variables $x=M+X$, $z=M+Z$
  and the elementary estimate
  $1+(M+Z)^2\le C(1+M^2)(1+Z^2)$, we have
  \begin{align*}
    \iint_{x,z\leq M}|V_{T,\lambda}(x,z)|^2 dx dz \leq C_T(1+M^2) e^{2T\lambda^2-2a^2\lambda}J_{\lambda, T},
  \end{align*}
  where
  \begin{align*}
    J_{\lambda, T}:= \iint_{X, Z\leq 0}(1+Z^2)\exp\left\{-\frac{(X+Z+2T\lambda)^2}{2T}
    \right\} dXd Z.
  \end{align*}
  Let $u=-(X+Z)$ and $s=-Z$ with $u\geq0$ and $0\leq s\leq u$. Then
  \begin{align*}
    J_{\lambda, T}&= \int_0^\infty e^{-\frac{(2T\lambda-u)^2}{2T}} \left(\int_0^u (1+s^2) ds\right)du\\
    &=\int_0^\infty (u+\frac{u^3}{3})e^{-\frac{(u-2T\lambda)^2}{2T}} du = \int_{-2T\lambda}^\infty (v+2T\lambda+\frac{(v+2T\lambda)^3}{3})e^{-\frac{v^2}{2T}} dv\\
    &\leq C_T \int_\R (1+|v|^3+\lambda^3) e^{-\frac{v^2}{2T}}dv.
  \end{align*}
  Thus, it follows that
  \begin{equation}\label{eq:V-ref-HS}
    \iint_{x,z\le M}|V_{T,\lambda}(x,z)|^2\,dx\,dz
    \le C_T(1+M)^2(1+\lambda^3)
    e^{2T\lambda^2-2a^2\lambda}.
  \end{equation}

  On the complementary region at least one of $X,Z$ is nonnegative and $Q_{M,T}=p_T(x-z)$ by~\eqref{eq:Q-piecewise-L21}.  We consider the integral $\int_0^\infty d X \int_\R d Z |V_{T,\lambda}(X+M, Z+M)|^2$ first.  Using the identity
  \begin{align*}
    -\frac{(X-Z)^2}{2T}-2\lambda(X+Z)-4\lambda M
    &=-\frac{(Z-X+2T\lambda)^2}{2T}
    +2T\lambda^2-4\lambda X-4\lambda M
  \end{align*}
  we write
  \begin{align*}
    & \int_0^\infty d X \int_\R d Z |V_{T,\lambda}(X+M, Z+M)|^2 \\
    &\qquad \qquad = \frac{e^{2T\lambda^2-4\lambda M}}{4\pi T}
    \int_0^\infty d X \int_\R d Z (1+(Z+M)^2)e^{-\frac{(Z-X+2T\lambda)^2}{2T}
    -4\lambda X}\\
    &\qquad \qquad = \frac{e^{2T\lambda^2-4\lambda M}}{4\pi T}
    \int_0^\infty d X \int_\R d Z (1+(Z+X-2T\lambda+M)^2)e^{-\frac{Z^2}{2T}
    -4\lambda X}\\
    &\qquad \qquad \leq C_Te^{2T\lambda^2-4\lambda M}  \int_0^\infty d X \int_\R d Z (1+Z^2+X^2+\lambda^2+M^2)e^{-\frac{Z^2}{2T}
    -4\lambda X}.
  \end{align*}
  Using the fact $\int_0^\infty e^{-4\lambda X}d X= c_1\lambda^{-1}$ and $\int_0^\infty X^2 e^{-4\lambda X}d X= c_2\lambda^{-3}$,  we obtain that for all $\lambda >T+1$,
  \begin{align*}
    \int_0^\infty d X \int_\R d Z |V_{T,\lambda}(X+M, Z+M)|^2  \leq C_T (1+M^2)(1+\lambda)e^{2T\lambda^2-4\lambda M}.
  \end{align*}
  Similarly, we also have
  \begin{align*}
    \int_0^\infty d Z \int_\R d X |V_{T,\lambda}(X+M, Z+M)|^2  \leq C_T (1+M^2)(1+\lambda)e^{2T\lambda^2-4\lambda M}.
  \end{align*}
  As a consequence, we obtain
  \begin{equation*}
    \iint_{x>M \, \, \text{or}\, \, z>M}|V_{T,\lambda}(x,z)|^2\,dx\,dz
    \le C_T (1+M^2)(1+\lambda)e^{2T\lambda^2-4\lambda M}
  \end{equation*}
  which together with~\eqref{eq:V-ref-HS} proves Lemma~\ref{lem:V}.
  This completes the proof of Lemma~\ref{lem:V}.
\end{proof}

\subsection{Weighted heat operators}
\begin{lemma}\label{lem:2to2}
  Let $U_a$ be defined in~\eqref{eq:UPhi} for $a>0$. Then we have for $L>0$
  \begin{equation*}\label{eq:102}
    \norm{U_a P_L U_a^{-1}}_{2\to2}\leq e^{La^2}.
  \end{equation*}
\end{lemma}
\begin{proof}
  We write
  \begin{align*}
    (U_aP_LU_a^{-1}f)(x)
    &=e^{-ax}\int_{\R}p_L(x-y)e^{ay}f(y)\,dy
    =\int_{\R}e^{-a(x-y)}p_L(x-y)f(y)\,dy.
  \end{align*}
  Thus, $U_aP_LU_a^{-1}$ is the convolution operator
  $
  U_aP_LU_a^{-1}f=k_{L,a}*f,
  $
  with
  \begin{equation*}\label{eq:kernel}
    k_{L,a}(v)=e^{-av}p_L(v).
  \end{equation*}
  By Young's convolution inequality,
  $
  \norm{k_{L,a}*f}_{2}\le \norm{k_{L,a}}_1\norm{f}_2.
  $
  Hence, \begin{equation*}\label{eq:young}
    \norm{U_aP_LU_a^{-1}}_{2\to2}\le \norm{k_{L,a}}_1.
  \end{equation*}
  We now compute the $L^1$ norm of $k_{L, a}$ exactly. Since $p_L\ge0$,
  \begin{align*}
    \norm{k_{L,a}}_1
    &=\int_{\R}e^{-av}p_L(v)\,dv
    =\frac{1}{\sqrt{4\pi L}}\int_{\R}
    \exp\!\left(-\frac{v^2}{4L}-av\right)\,dv.
  \end{align*}
  Complete the square:
  $
  -\frac{v^2}{4L}-av
  =-\frac{(v+2La)^2}{4L}+La^2.
  $
  Therefore, \begin{align*}
    \norm{k_{L,a}}_1
    &=e^{La^2}\frac{1}{\sqrt{4\pi L}}
    \int_{\R}\exp\!\left(-\frac{(v+2La)^2}{4L}\right)\,dv
    =e^{La^2}.
  \end{align*}
  This completes the proof of Lemma~\ref{lem:2to2}.
\end{proof}

\begin{lemma}\label{lem:UPhi}
  For each $T>0$ there are constants $C_T,K>0$ such that, for all $M,L>0$ with $a=\sqrt{2M}\ge T+1$, and $U_a,\Phi$ as in~\eqref{eq:UPhi},
  \begin{align*}
    \|U_aQ_{M, T}P_LU_a\Phi\|_2 \leq C_T(1+M+L)^K
    \exp\{La^2+Ta^2-a^3\}.
  \end{align*}
\end{lemma}
\begin{proof}
  We write
  \begin{align}\label{eq:UB}
    U_aQ_{M, T}P_LU_a\Phi
    =(U_aQ_{M,T}U_a\Phi)\,\mathcal B_{L,a}, \qquad \text{with}\,\,\, \mathcal B_{L,a}:=\Phi^{-1}U_a^{-1}P_LU_a\Phi.
  \end{align}
  The kernel of $\mathcal B_{L,a}$  is given by
  \[
    \mathcal B_{L,a}(z,w)
    =\frac{\phi(w)}{\phi(z)}e^{a(z-w)}p_L(z-w).
  \]
  Since $\phi$ is $1$-Lipschitz and $\phi\ge1$, the elementary inequality
  $
  \frac{\phi(w)}{\phi(z)}\le 1+|z-w|
  $
  holds for all $z,w\in\R$.  Define
  \[
    k(v):=(1+|v|)e^{av}p_L(v).
  \]
  Then $|\mathcal B_{L,a}(z,w)|\le k(z-w)$, and therefore, for
  $f\in L^2(\R)$,
  \[
    |\mathcal B_{L,a}f(z)|\le (k*|f|)(z).
  \]
  Young's convolution inequality gives
  \[
    \|\mathcal B_{L,a}f\|_2
    \le \|k*|f|\|_2
    \le \|k\|_1\|f\|_2,
  \]
  which implies
  \begin{align*}
    \|\mathcal B_{L,a}\|_{2\to2}
    &\le \|k\|_1
    =\int_{\R}(1+|v|)e^{av}p_L(v)\,dv.
  \end{align*}
  Using the equality
  $
  -\frac{v^2}{4L}+av
  =-\frac{(v-2La)^2}{4L}+La^2,
  $
  we see that $e^{av}p_L(v)$ is $e^{La^2}$ times the density of
  $G\sim {\rm N}(2La,2L)$, and hence
  \[
    \int_{\R}(1+|v|)e^{av}p_L(v)\,dv
    =e^{La^2}(1+\mathbb E|G|)
    \le C(1+aL+\sqrt L)e^{La^2}.
  \]
  Thus, we deduce that
  \begin{align}
    \|\mathcal B_{L,a}\|_{2\to2}
    &\le C(1+aL+\sqrt L)e^{La^2}.
    \label{eq:B-Schur-final}
  \end{align}
  By Lemma~\ref{lem:V}, we conclude from~\eqref{eq:B-Schur-final} and~\eqref{eq:UB} that
  \begin{equation*}\label{eq:QPL-final}
    \|U_aQP_LU_a\Phi\|_2
    \le
    C_T(1+M+L)^K
    \exp\{La^2+Ta^2-a^3\}.
  \end{equation*}
  This completes the proof of Lemma~\ref{lem:UPhi}.
\end{proof}

\subsection{A truncated Gaussian convolution}
\begin{lemma}\label{lem:J}
  Fix $M, L>0$. Then for $x\geq M$ and $w\leq M$, we have
  \begin{equation*}\label{eq:FR-J-exact}
    \int_{-\infty}^M p_L(x-z)p_1(2M-z-w)\,dz
    =p_{L+1}(\alpha-\beta)
    \Phi_{\rm G}\!\left(-\frac{\alpha+L\beta}{\sqrt{2L(L+1)}}\right),
  \end{equation*}
  where $\Phi_{\rm G}$ is the standard normal distribution function and $\alpha=x-M$, $\beta=M-w$.
\end{lemma}
\begin{proof}
  Denote $S=L+1$ and write
  $
  x=M+\alpha, w=M-\beta,
  \alpha,\beta\ge0.
  $
  For fixed $\alpha,\beta$,  define
  \[
    J(\alpha,\beta)
    :=\int_{-\infty}^M p_L(x-z)p_1(2M-z-w)\,dz.
  \]
  We now compute this truncated convolution explicitly.  Put
  $
  c:=2M-w
  $ and write
  \[
    p_L(x-z)p_1(2M-w-z)
    =\frac{1}{4\pi\sqrt L}
    \exp\left\{-\frac{(z-x)^2}{4L}
    -\frac{(z-c)^2}{4}\right\}.
  \]
  Because
  $S=L+1$, we have
  \begin{align*}
    \frac{(z-x)^2}{L}+(z-c)^2
    &=\frac{S}{L}\left(z-\frac{x+Lc}{S}\right)^2
    +\frac{(x-c)^2}{S}.
  \end{align*}
  Thus, with
  $
  m:=\frac{x+Lc}{S}=\frac{x+L(2M-w)}{S},
  $
  we have the pointwise factorization
  \begin{align*}
    p_L(x-z)p_1(c-z)
    &=\frac{1}{4\pi\sqrt L}
    \exp\left\{-\frac{(x-c)^2}{4S}\right\}
    \exp\left\{-\frac{S}{4L}(z-m)^2\right\}\\
    &=p_S(x-c)\,
    \sqrt{\frac{S}{4\pi L}}
    \exp\left\{-\frac{S}{4L}(z-m)^2\right\}.
  \end{align*}
  Integrating over $z\le M$ therefore gives
  \begin{align*}
    J(\alpha,\beta)
    &=p_S(x-c)
    \int_{-\infty}^{M}
    \sqrt{\frac{S}{4\pi L}}
    \exp\left\{-\frac{S}{4L}(z-m)^2\right\}\,dz \\
    &=p_S(x-c)\,
    \Phi_{\rm G}\left(\frac{M-m}{\sqrt{2L/S}}\right).
  \end{align*}
  Here $\Phi_{\rm G}$ is the standard normal distribution function.  We next rewrite
  both quantities in terms of $\alpha,\beta$.  Since
  $
  x=M+\alpha$ and $c=2M-w=M+\beta,
  $
  we have
  $
  x-c=\alpha-\beta,
  $
  so that the first factor is $p_S(\alpha-\beta)$.  Also
  \begin{align*}
    m-M
    &=\frac{x+Lc-MS}{S}
    =\frac{M+\alpha+L(M+\beta)-M(L+1)}{S}
    =\frac{\alpha+L\beta}{S}.
  \end{align*}
  Hence, $
  \frac{M-m}{\sqrt{2L/S}}
  =-\frac{\alpha+L\beta}{\sqrt{2LS}}.
  $ This proves Lemma~\ref{lem:J}.
  This completes the proof of Lemma~\ref{lem:J}.
\end{proof}

\noindent
\textbf{Acknowledgments}.
L.C. was supported in part by NSF grants DMS-2246850/2443823 and a Collaboration
Grant for Mathematicians (\#959981) from the Simons Foundation. F.P. was
supported in part by National Key R\&D Program of China (No. 2022YFA 1006500) and National Natural Science Foundation of China (No. 12571153).  The second author would like to thank Xinxin Chen and Hui He for raising the question of fluctuation of the maximum of the Airy$_1$ process during a presentation on the work~\cite{Pu2025}. We acknowledge the use of ChatGPT.

\bigskip
\noindent
\textsc{Le Chen, Department of Mathematics and Statistics, Auburn University, Auburn, AL 36849, USA.}

\noindent
\textit{Email address:} \texttt{lzc0090@auburn.edu}

\medskip
\noindent
\textsc{Fei Pu, Laboratory of Mathematics and Complex Systems,
School of Mathematical Sciences, Beijing Normal University, 100875, Beijing, China}

\noindent
\textit{Email address:} \texttt{fei.pu@bnu.edu.cn}

\end{document}